\documentclass[10pt,a4paper,leqno]{amsart}
\usepackage{amssymb}
\usepackage{graphicx}
\usepackage{pdfsync}
\usepackage{xspace}
\usepackage{fancyhdr}
\usepackage[english]{babel}
\usepackage{amsfonts}
\usepackage{latexsym}
\usepackage{graphics}
\usepackage{hyperref}
\usepackage{pifont,bbding}
\usepackage{color}
\usepackage{bm}
\usepackage{enumerate,amsthm,amsmath,comment,psfrag}
 \usepackage{mathrsfs} 
 \usepackage{nicefrac}

\newtheorem{remark}{Remark}[section]
\newtheorem{lemma}{Lemma}[section]
\newtheorem{proposition}{Proposition}[section]

\numberwithin{equation}{section}
\newcommand{\ep}{\varepsilon}

\newcommand{\ii}{\mathrm{i}}

\newcommand{\p}{\partial}

\newcommand{\xv}{\mathbf{x}}
\newcommand{\V}{\mathsf{v}}

\newcommand{\Dh}{\Delta_h}
\newcommand{\E}{\mathcal{E}}

\newcommand{\Pc}{\mathcal{P}}

\newcommand{\Rc}{\mathcal{R}e}

\newcommand{\Pb}{\mathbb{P}}
\newcommand{\Rb}{\mathbb{R}}

\newcommand{\Cb}{\mathbb{C}}

\newcommand{\Th}{\mathcal{T}_h}
\newcommand{\Vh}{\mathcal{V}_h}

\newcommand{\dif}{\,\text{d}}

\newcommand{\RED}[1]{{\color{black}{#1}}}

\makeatletter
\newcommand*\bigcdot{\mathpalette\bigcdot@{.5}}
\newcommand*\bigcdot@[2]{\mathbin{\vcenter{\hbox{\scalebox{#2}{$\m@th#1\bullet$}}}}}
\makeatother

\begin{document}
%
%
%-------------------------------------------------------------------------------
% Top matter
%-------------------------------------------------------------------------------
%
%
\title[A novel  relaxation scheme for Schr\"odinger-Poisson]
{
A novel, structure-preserving, second-order-in-time relaxation scheme for   Schr\"odinger-Poisson systems}
\date{\today}
\author{Agissilaos  Athanassoulis}
\address[Agissilaos  Athanassoulis]{Department of Mathematics, University of Dundee, Dundee DD1 4HN, Scotland, UK} \email{a.athanassoulis@dundee.ac.uk}

\author{Theodoros Katsaounis}
\address[Theodoros Katsaounis]{Dept. of Math. and  Applied Mathematics, Univ. of Crete, GREECE \& IACM--FORTH, Heraklion, GREECE
}
\email{thodoros.katsaounis@uoc.gr}
\author{Irene Kyza}
\address[Irene Kyza]{Department of Mathematics, University of Dundee, Dundee DD1 4HN, Scotland, UK} \email{ikyza@dundee.ac.uk}

\author{Stephen Metcalfe}
\address[Stephen Metcalfe]{Dept. of Mechanical Engineering, McGill Univ., Montreal, Canada} \email{smetcalfephd@gmail.com}

%\mu
\begin{abstract}
We introduce a new \RED{structure preserving}, second order in time relaxation-type scheme for approximating solutions of the Schr\"odinger-Poisson system. More specifically, we use the Crank-Nicolson scheme as a time stepping mechanism, whilst the nonlinearity is handled by means of a relaxation approach in the spirit of \cite{Besse, KK} for the nonlinear Schr\"odinger equation. For the spatial discretisation we use the standard conforming finite element scheme. The resulting scheme is explicit with respect to the nonlinearity, {\color{black} i.e. it requires the solution of a linear system for each time-step, and satisfies discrete versions  of the system's mass conservation and energy balance laws for constant meshes. The scheme is seen to be second order in time. We conclude by presenting some numerical experiments, including an example from cosmology and an example with variable time-steps} which demonstrate the effectiveness and robustness of the new scheme.
\end{abstract}
\maketitle

%\tableofcontents
%
%
%\nocite{*}
%-------------------------------------------------------------------------------
\section{Introduction}
%-------------------------------------------------------------------------------
%
\subsection{Statement of the problem}
%While many numerical methods exist for their simulation, challenges remain. In this paper, we introduce a second order in time, linearly implicit numerical scheme applicable to a broad class of  Schr\"odinger-Poisson type systems. This scheme satisfies discrete versions of the continuous mass conservation and energy balance laws. 

%\subsection{The continuous problem}
Schr\"odinger-Poisson-type systems appear in many applications, including  semiconductors \cite{RS,Karner,MRS}, plasma physics \cite{BTGIFG,Shukla}, optics \cite{Paredes} and cosmology \cite{UKH,KVS,WK,DW}. 
In this paper we consider a  class of Schr\"odinger-Poisson systems (SPS), namely the following initial-boundary value problem, either with homogeneous Dirichlet boundary conditions %(case $\mu=0$ in system \eqref{GeSP}),
or periodic boundary conditions: %(case $\mu=\|u_0\|_{L^2}^2$ in system \eqref{GeSP}). That is, given an initial time  $\tau\ge 0$ and a final time $T>\tau,$ 
We seek  a \emph{wavefunction} $u : \varOmega \!\times \!(\tau,T) \to \Cb$ and the associated \emph{potential} $\V : \varOmega \!\times \!(\tau,T) \to \Rb$ such that
\begin{equation}
\label{GeSP} \left \{
\begin{aligned}
&u_t-\ii p(t)\Delta u +\ii q(t) \V u=0,    &&\quad\mbox{in ${\varOmega}\!\times\! (\tau,T)$,}& \\
&\Delta\V = |u|^2 - \mu, &&\quad\mbox{in ${\varOmega}\!\times\! (\tau,T)$,}&\\
&u(\xv,{\color{black}\tau})=u_0(\xv), &&\quad\mbox{in ${\varOmega}$}, & \\
&\{ \mu = 0 \text{ and }  u = \V = 0\}, \ \text{ OR } \ \{\mu = \|u_0\|_{L^2}^2 \text{ and }  u, \V \text{ periodic,\} }&&\quad\mbox{on $\partial \varOmega\!\times\! (\tau,T]$,}&  \\
%&\mu = \|u_0\|_{L^2}^2 :  u, \V \text{ periodic, } &&\quad\mbox{on $\partial \varOmega\!\times\! (\tau,T]$.}&  \\
\end{aligned}
\right.
\end{equation}%

\noindent
In \eqref{GeSP}, $\tau\ge 0$ is a given initial time and $T>\tau $ is a given final time. %For $\mu=0$, we get homogeneous Dirichlet boundary conditions, while for  $\mu=\|u_0\|_{L^2}^2$ we have periodic boundary conditions.
 The domain $\varOmega \subset \Rb^d, \ d=1,2,3$, is assumed to be bounded, convex and polygonal in the Dirichlet case, and a $d-$dimensional parallelepiped in the periodic case.  The normalisation of $\mu$ ensures that the elliptic problem is well-posed in both cases.  For the initial condition we have  $u_0 \in H^1(\varOmega)$. The  coefficients $p(t), \ q(t)$ are smooth and real valued; the introduction of time-dependent coefficients is directly motivated by the cosmological application \cite{UKH,KVS}. More details and simulations of that problem can be found in Section \ref{sec:cosmoex}.
 %The final time $T \in \mathbb{R}^+$ while the data is assumed to satisfy $\alpha, \varepsilon \in \mathbb{R}^+$, $\beta \in \mathbb{R}$ and $u_0 \in H^1(\varOmega)$.
 
System \eqref{GeSP} satisfies mass conservation and energy balance laws that are of great physical relevance. These are discussed in detail in Section \ref{sec:energybalance}. {\color{black}Proposing schemes that satisfy discrete analogues of these laws is a significant goal, as typically this leads to good qualitative behaviour of numerical solutions for longer computational times. Moreover, when these discrete laws are verified unconditionally in the time-step size, this provides flexibility for dealing with stiffness issues. Structure preserving schemes for \eqref{GeSP} are of great physical relevance since in recent years the SPS system is used extensively as an alternative to the computationally expensive Vlassov-Poisson system with applications in cosmology, see e.g. \cite{KVS, UKH, WK} and the references therein.
 
 Our goal in this paper is to propose a scheme that is is linearly implicit, unconditionally structure preserving (in the sense of satisfying discrete energy and mass balance laws without restriction on the time-step size) and second order accurate in time. An important advantage of linearly implicit schemes is that they are faster. Moreover, they are easier to implement, as no iterative scheme (Newton, fixed point etc) has to be selected and calibrated. Finally, by being more directly implementable, linearly implicit schemes are more amenable to further work on a posteriori error control.}

%\subsection{More on the Schr\"odinger-Poisson system}

An important special case of \eqref{GeSP} is the problem with constant coefficients and homogeneous Dirichlet boundary conditions; 
 taking $ p(t) = \dfrac{\ep}{2\alpha^2}$, $ q(t) = \dfrac{\beta}{\ep\alpha}$ and $\mu=0$  we obtain the following standard Schr\"odinger-Poisson system, 
 \begin{equation}
\label{SP} \left \{
\begin{aligned}
&u_t-\frac{i\ep}{2\alpha^2}\Delta u +\frac{i\beta}{\ep\alpha} \V u=0,    &&\quad\mbox{in ${\varOmega}\!\times\! (\tau,T)$,}& \\
&\Delta\V =  |u|^2,   &&\quad\mbox{in ${\varOmega}\!\times\! (\tau,T)$,}&\\
&u=0,\ \V=0,  &&\quad\mbox{on $\partial \varOmega\!\times\! (\tau,T]$,}&  \\
&u(\xv,{\color{black}\tau})=u_0(\xv), &&\quad\mbox{in ${\varOmega}$},& \\
\end{aligned}
\right.
\end{equation}%
where the  parameter $\ep$ represents the ratio of the Planck constant to the mass of the particle, while $\alpha>0$ and $\beta\in\Rb$ are given constants. When $\ep$ is small, this is called the semiclassically scaled  problem, and 
 it is expected formally that as $\varepsilon\to 0^+$  the Schr\"odinger-Poisson system \eqref{SP} approximates, in some sense,  the classical Vlasov-Poisson equations, cf., e.g., \cite{ZZM}. 
 
 Concerning the existence and uniqueness of solutions to the Schr\"odinger-Poisson system \eqref{SP}, most of the analytical results are for the full space case $\varOmega=\mathbb{R}^d;$ we refer to the works \cite{BM,Castella,IZL} and the book \cite{Cazenave}.  In \cite{AMP}, the authors analyze a  transient Schr\"odinger-Poisson system with transparent boundary conditions while in \cite{TM} the stationary spherically symmetric case was analyzed.  The asymptotic behaviour of solutions to the  Schr\"odinger-Poisson system is studied in \cite{AS} via a variational approach. 
%

%\noindent In the field of cosmology, system \eqref{SP} can be considered as an alternative  model for  describing collisionless self-gravitating matter, \cite{KVS}. In this case,  $\mu = 1$ and the parameter $\alpha$ depends on time resulting in  a system of the form \eqref{GeSP} with $p(t) = \frac{\ep}{2 t^{\nicefrac{3}{2}}}$, $q(t) = \frac{\beta}{\ep t^{\nicefrac{1}{2}}}$. 

\subsection{Existing numerical methods for the Schr\"odinger-Poisson system}
{\color{black} There exists a very large literature for the numerical approximation of the nonlinear Schr\"odinger equation with power nonlinearity (NLS),  cf. e.g. \cite{Besse,Besse2,KK,AkrivisLi,ADK,ADKM,Fei,Delfour,Bao,Kara,KM2,Berland,Chartier,Deg1,Deg2,Deg3,Besse3,Hederi,Tha,Zouraris1} and the references therein for a sample of such works. In contrast, the numerical methods available for the SPS are not that many. In what follows, we focus on  numerical methods for the SPS \eqref{GeSP}. } 
%\BLUE{, and do not comment in detail on the NLS literature that does not directly apply to the SPS \eqref{GeSP} or to some simpler version of it like \eqref{SP}.}

%For the SPS, due to the nonlinear interactions and the scaling $\varepsilon\ll 1 $ (cf. \eqref{SP}), the interesting physical quantities of \eqref{GeSP}  develop sharply localised features, making the proposition of such a method a non-trivial task. 
\RED{The existing methods in the literature for SPS are with uniform temporal and mesh sizes and/or the nonlinear term is treated implicitly. This leads to practical difficulties, which might  be the reason that in the literature most of the numerical implementations for the SPS \eqref{GeSP} are performed in the one-dimensional case. 

}

%In the literature, there are several methods available for the numerical approximation of solutions to the Schr\"odinger-Poisson system \eqref{SP}. 

More precisely, popular methods in the literature for the approximation of the SPS \eqref{GeSP}, include the Crank-Nicolson, the Gaussian beam or the time-splitting method for the time-discretization, while finite differences or spectral methods are used for the spatial discretization. In \cite{AKKT}, the authors conduct an error and stability analysis for an operator splitting finite element discretization of \eqref{SP} whilst an error analysis for the semidiscrete Galerkin finite element scheme is presented in \cite{BILZ}. Utilizing a Crank-Nicolson temporal and finite difference spatial discretization of \eqref{SP}, a predictor-corrector scheme is studied in \cite{RS} and the spherically symmetric case is studied in \cite{EZ}. In \cite{BNS}, the behaviour of the solution of the Schr\"odinger-Poisson-X$\alpha$ system is explored through a discretization based on the time splitting spectral method. A time semidiscrete scheme for \eqref{SP} using Strang splitting is studied extensively in \cite{Lubich} and an error analysis is provided. The Gaussian beams method is introduced in \cite{JWY} for the numerical simulation of \eqref{SP} in the one dimensional case whilst in \cite{Zhang} error estimates are obtained for a Crank-Nicolson in time, compact finite difference in space discretization of \eqref{SP}. A numerical method consisting of a backward Euler in time, pseudo-spectral method in space is utilized in \cite{ZD} to approximate the ground states and the solution of the Schr\"odinger-Poisson-Slater system (which also includes \eqref{SP}). A spectral discontinuous Galerkin method in space coupled with a Runge-Kutta scheme in time is used to study solutions of \eqref{SP} in \cite{LC}. 

{\color{black} All of the  aforementioned methods for the SPS are implicit with respect to the nonlinear term, which means that a nonlinear system must be solved for their implementation. Also none of the above methods satisfy a discrete analogue of the energy balance (or the energy conservation) that the continuous problem does. \emph{To the best of our knowledge, at the moment, the literature lacks  a numerical method for the Schr\"odinger-Poisson system that  satisfies a discrete analogue of the energy balance of the continuous problem.}}

\subsection{Main contribution of this paper}
\emph{In this paper we introduce a novel numerical scheme for the discretization of the SPS \eqref{GeSP}.} In particular, for the spatial discretization we use the standard conforming finite element method, while for the temporal discretization we propose a new Crank-Nicolson relaxation-type method. More precisely, we linearize the nonlinear term $|u|^2-\mu$ in the potential equation of \eqref{GeSP} and discretize it using the Crank-Nicolson method. The main advantage of the particular scheme is that it avoids solving a computationally expensive nonlinear system, whilst maintaining the second order temporal accuracy of the Crank-Nicolson method. For the proposition of this method we were inspired by Besse, \cite{Besse,Besse2}, who used a relaxation-type Crank-Nicolson  scheme with constant time-steps for the nonlinear Schr\"odinger equation. Katsaounis and Kyza in \cite{KK} generalised the scheme of \cite{Besse} to variable time-steps and used it to derive a posteriori error estimates, while Zouraris in \cite{Zouraris} used a similar scheme for a semilinear parabolic equation.

In particular, \emph{our main contribution  is the proposition of a fully discrete relaxation-type Crank-Nicolson  finite element scheme with variable temporal and spatial mesh sizes, which is implicit with respect to the Laplacian terms of \eqref{GeSP} (hence guaranteeing stability), but explicit in the nonlinearity of the potential equation of \eqref{GeSP}}. The fact that the proposed method is with variable temporal and mesh-sizes opens the road for the a posteriori error analysis and the proposition of adaptive algorithms for \eqref{GeSP}; this is the focus of a forthcoming paper. \emph{Our method is  easily implementable and is numerically shown to be of second order of accuracy. Moreover, it is  proven to inherit on the discrete level the mass conservation and energy balance laws of the original continuous problem, for constant meshes. To the best of our knowledge, this is the first time that such a method, satisfying all the above mentioned properties, is proposed in the literature for the SPS \eqref{GeSP}. {\color{black}In fact, to the best of our knowledge this is the first method proposed for the SPS that satisfies a discrete energy law.}}

{\color{black}
\subsection{Idea behind the scheme}
The scheme is presented in full detail in Section 3, however at this point we can give a heuristic description of the main idea behind it.  For example, a semi-discrete fully implicit Crank-Nicolson scheme with constant time-step $k$ for the standard SPS \eqref{SP} could be written as
 \begin{equation}
\label{SPmid} \left \{
\begin{aligned}
&\frac{U^{n}-U^{n-1}}{k}-\frac{i\ep}{2\alpha^2}\Delta U^{n-\frac{1}2} +\frac{i\beta}{\ep\alpha} V^{n-\frac{1}2} U^{n-\frac{1}2}=0,    && & \\
&\Delta V^{n-\frac{1}2} =  |U^{n-\frac{1}2}|^2,   && &
\end{aligned}
\right.
\end{equation}%
with the usual notation of $U^{n-\frac{1}2}:=(U^{n}+U^{n-1})/2,$ $V^{n-\frac{1}2}:=(V^{n}+V^{n-1})/2.$ This might be an interesting scheme (second order in $k$, mass preserving), but each time-step involves a nonlinear system in the unknowns $U^n,$ $\V^n$ and would be computationally expensive and complicated to implement. 

So instead of that, the idea here is to use an auxiliary variable on a staggered timegrid. We introduce an auxiliary variable $\Phi^{n-\frac{1}2}$ as a proxy for the position density and update it by linear extrapolation
\begin{equation}\label{eq:interp111}
\frac{\Phi^{n-\frac{1}2} + \Phi^{n-\frac{3}2}}2 = |U^{n-1}|^2,
\end{equation}
then define $V^{n-\frac{1}2}$ simply as
\[
\Delta V^{n-\frac{1}2} = \Phi^{n-\frac{1}2}.
\]
Thus we end up with a semi-discrete scheme of the form
 \begin{equation}
\label{SPmid} \left \{
\begin{aligned}
&\Phi^{n-\frac{1}2}  = 2|U^{n-1}|^2 - \Phi^{n-\frac{3}2}, &&& \\
&\Delta V^{n-\frac{1}2} =  \Phi^{n-\frac{1}2},   && &\\
&\frac{U^{n}-U^{n-1}}{k}-\frac{i\ep}{2\alpha^2}\Delta U^{n-\frac{1}2} +\frac{i\beta}{\ep\alpha} V^{n-\frac{1}2} U^{n-\frac{1}2}=0,   && & 
\end{aligned}
\right.
\end{equation}%
where the equations above appear in the order they would be solved. 
For variable time-step, the same idea yields the scheme in \eqref{CNrelax} by adjusting the linear extrapolation step of \eqref{eq:interp111}.
}
%In this work, we introduce a novel numerical scheme for the discretization of the Schr\"odinger-Poisson system \eqref{SP}. We use the standard continuous finite element method for the spatial discretization and the Crank-Nicolson method for the time stepping mechanism. To handle the nonlinearity, we utilize the relaxation scheme  introduced by Besse \cite{Besse} for the nonlinear Schr\"odinger equation which offers a way of  linearizing the nonlinear term in the potential equation \eqref{SP}. For the Besse discretization of the nonlinear Schr\"odinger equation, the authors derived a-posteriori error estimates up to the critical exponent in \cite{KK}.  The main advantage of the scheme is that we avoid solving a computationally expensive nonlinear system whilst maintaining the second order temporal accuracy of the Crank-Nicolson method. The proposed method conserves the density at the discrete level and a energy balance is satisfied as well. 
\bigskip

\subsection{Organization of the paper}
The rest of the paper is organized as follows: In Section 2 we derive the balance laws of the Schr\"odinger-Poisson system \eqref{GeSP} while in Section 3 we introduce our new numerical method and derive discrete variants of the system's conservation laws. Section 4 explains the practical implementation of the numerical scheme. In Section 5, we present numerical experiments which verify the accuracy and efficiency of the method. We finish the section by applying the numerical  method to a concrete example from cosmology. %Finally, we draw conclusions in Section 6.

%------------------------------------------------------------------------------- 
\section{Mass conservation \& Energy balance}\label{sec:energybalance}
%-------------------------------------------------------------------------------
%
The standard Schr\"odinger-Poisson system \eqref{SP} exhibits mass and energy conservation, consistently with its quantum mechanical interpretation \cite{BEGMY}. If however the coefficients $p(t),q(t)$ vary in time, the more general system \eqref{GeSP} \RED{still satisfies mass conservation} and exhibits a precise {\em energy balance law} instead of simple conservation.
More specifically, given \eqref{GeSP} let us define the \emph{mass} 
\begin{align}
\mathcal{M}(t) & := \|u(t)\|^2, \label{cldef1}
\end{align}
and the 
 \emph{kinetic energy} $\E_k(t)$ and \emph{potential energy} $ \E_{\V}(t)$  
\begin{align}
\E_k(t) & := \|\nabla u(t)\|^2, \quad  \E_{\V}(t) := \|\nabla \V(t)\|^2 \ = - \int_{\varOmega} \V(t) (|u(t)|^2-\mu) \dif x, \label{cldef2} 
\end{align}
where    $\|\cdot\|$ denotes the $L^2$-norm over $\varOmega$.  
It is worth noting that  the integration by parts, %e.g. such as
\begin{equation}
\int_\varOmega \nabla \V(t) \cdot \nabla \V(t) \dif x = - \int_{\varOmega} \V(t) \Delta \V(t) \dif x,
\end{equation}
which is used in showing that the two expressions for the potential energy  in \eqref{cldef2} are equivelant, works equally well with both the periodic or the homogeneous Dirichlet boundary conditions. With this at hand, we are now ready to prove the following:
\begin{lemma}[Continuous Mass Conservation \& Energy Balance]\label{conslaws} If $(u, \V)$ is a solution of  \eqref{GeSP}  then , for $ \tau \le  t \le T$ 
\begin{align}
& \mathcal{M}(t)  = \mathcal{M}(\tau), &&\text{Conservation of Mass}, \label{cl1}\\
& p(t) \frac{\dif}{\dif t} \E_k(t)  - \frac{q(t)}{2} \frac{\dif}{\dif t} \E_{\V}(t)=0,  &&\text{Balance of Energy}. \label{cl2}
\end{align}
\end{lemma}
\begin{proof} We begin the proof by deriving conservation of mass. To obtain this, we multiply the Schr{\"o}dinger equation by $\bar{u}$ and integrate over $\varOmega$ yielding
\begin{equation}
\begin{aligned}
\int_{\varOmega} \bar{u}u_t  \dif x +\ii p(t) \|\nabla u\|^2 + \ii q(t)\int_{\varOmega} \V|u|^2 \dif x &= 0.
\end{aligned}
\end{equation}
Since $p(t), \ q(t), \ \V$ are real-valued, taking real parts immediately implies that 
\begin{equation}
\begin{aligned}
\int_{\varOmega} (\bar{u}u_t + u \bar{u}_t) \dif x=0 \,\,\implies \,\,
\frac{\dif }{\dif t}\mathcal{M}(t) = 0.
\end{aligned}
\end{equation}
and \eqref{cl1} readily holds.
%Hence, upon integration over $(t_0,t)$, we obtain the equation for conservation of mass:
%\begin{equation}
%\begin{aligned}
%\mathcal{D}(t) = \mathcal{D}(t_0).
%\end{aligned}
%\end{equation}

To derive the energy balance, we begin by differentiating the potential equation with respect to $t$ yielding
\begin{equation}
\begin{aligned}
\Delta \V_t = \frac{\partial}{\partial t}\bigg(|u|^2 - \mu \bigg) = 2\Rc(u\bar{u}_t).
\end{aligned}
\end{equation}
Multiplying this by $\V$ and integrating over $\varOmega$ we get
\begin{equation}
\begin{aligned}\label{encon1}
-\frac{1}{2}\frac{\dif}{\dif t}\|\nabla \V\|^2 = \int_{\varOmega} \V\Delta \V_t \dif x = 2 \!\int_{\varOmega} \V\Rc(u\bar{u}_t) \,\dif x.
\end{aligned}
\end{equation}
If we now resort to the Schr{\"o}dinger equation, multiply it by $\bar{u}_t$ and integrate over $\varOmega$ we obtain
\begin{equation}
\begin{aligned}
\int_{\varOmega} u_t \bar{u}_t \dif x - \ii p(t)\int_{\varOmega} \bar{u}_t \Delta u \dif x + \ii q(t) \int_{\varOmega} \V u\bar{u}_t \dif x = 0,
\end{aligned}
\end{equation}
or equivalently,
\begin{equation}
\begin{aligned}
\|u_t\|^2 +\ii p(t)\int_{\varOmega}\nabla u\bigcdot \nabla \bar{u}_t \dif x + \ii q(t) \int_{\varOmega} \V u\bar{u}_t \dif x = 0.
\end{aligned}
\end{equation}
Taking imaginary parts yields
\begin{equation}
\begin{aligned}
\frac12 p(t)\frac{\dif}{\dif t}\|\nabla u\|^2 + q(t) \int_{\varOmega} \V\Rc(u\bar{u}_t) \dif x = 0.
\end{aligned}
\end{equation}
We now substitute in \eqref{encon1} to get
\begin{equation}
\begin{aligned}\label{enconeq1}
% \frac{\dif \E}{\dif t} = \frac{\dif}{\dif t}\bigg(\frac{\ep^2}{\alpha^2}||\nabla u||^2 - \frac{\alpha}{\beta} ||\nabla \V||^2 \bigg)= 0.
p(t)\frac{\dif}{\dif t}\|\nabla u\|^2 - \frac12 q(t) \frac{\dif}{\dif t} \|\nabla \V\|^2 = 0.
\end{aligned}
\end{equation}
Using the second expression for the potential energy, namely 
\begin{equation}
\begin{aligned}
||\nabla \V||^2 = - \int_{\varOmega} \V \Delta \V \dif x = - \int_{\varOmega} \V ( |u|^2 - \mu) \dif x,
\end{aligned}
\end{equation}
we can equivalently get
\begin{equation}
\begin{aligned}\label{enconeq2}
%\frac{\dif \E}{\dif t} = \frac{\dif}{\dif t}\bigg(\frac{\ep^2}{\alpha^2}||\nabla u||^2 + \int_{\varOmega} \V |u|^2 \dif x \bigg) = 0.
p(t)\frac{\dif}{\dif t}\|\nabla u\|^2 + \frac12 q(t) \frac{\dif}{\dif t} \int_{\varOmega} \V (|u|^2-\mu) \dif x = 0.
\end{aligned}
\end{equation}
From either \eqref{enconeq1} or \eqref{enconeq2} we obtain the balance of energy \eqref{cl2}.
\end{proof}

\begin{remark}
From the previous lemma it's obvious that for $p,q$ constant in time one recovers the energy conservation law
\[
\frac{\dif }{\dif t} \left( p \mathcal{E}_k(t) {\color{black}-} \frac{q}2 \mathcal{E}_\V (t)\right)=0 \,\, \implies \,\, p \mathcal{E}_k(t) {\color{black}-} \frac{q}2 \mathcal{E}_\V (t)=p \mathcal{E}_k(\tau) {\color{black}-} \frac{q}2 \mathcal{E}_\V (\tau), \quad \tau\le t\le T.
\]
\end{remark}
%---------------------------------------------------------------------------------------
\section{A new Relaxation-type Numerical Method \& Discrete Balance Laws}\label{timediscr}
%--------------------------------------------------------------------------------------
%
The numerical scheme we propose here is inspired by the Crank-Nicolson relaxation method introduced by Besse in \cite{Besse} for the nonlinear Schr\"odinger equation. The first stage in the creation of a Besse-style relaxation scheme is to rewrite system \eqref{GeSP} via the introduction of an auxiliary variable $\phi = |u|^2-\mu$ which takes the place of the nonlinearity. We are thus now searching for a solution $(u, \V, \phi)$ of the following enlarged Schr\"odinger-Poisson system
\begin{equation}
\label{SPRLX} \left \{
\begin{aligned}
&u_t-\ii p(t)\Delta u +\ii q(t)\V u=0,    &&\qquad\mbox{in ${\varOmega}\!\times\! (\tau,T)$,}& \\
&\Delta\V = \phi,   &&\qquad\mbox{in ${\varOmega}\!\times\! (\tau,T)$,}&\\
&\phi = |u|^2 - \mu,  &&\qquad\mbox{in ${\varOmega}\!\times\! (\tau,T)$,}& \\
%&u=0,\ \V=0, \ \phi = 0, &&\qquad\mbox{on $\partial \varOmega\!\times\! (0,T]$,}&  \\
&u(\xv,{\color{black}\tau})=u_0(\xv), &&\qquad\mbox{in ${\varOmega}$},& \\
&\{ \mu = 0 \text{ and }  u = \V = 0\}, \ \text{ OR } \ \{\mu = \|u_0\|_{L^2}^2 \text{ and }  u, \V \text{ periodic, } \}&&\qquad\mbox{on $\partial \varOmega\!\times\! (\tau,T]$,}&  \\
%
%&\mu = \|u_0\|^2 :  u, \V \text{ periodic, }  &&\quad\mbox{on $\partial \varOmega\!\times\! (\tau,T]$,}&  \\
\end{aligned}
\right.
\end{equation}
which is, obviously, equivalent to the original problem \eqref{GeSP}. The numerical scheme that we will introduce in the sequel is based upon this enlarged formulation of the Schr\"odinger-Poisson system. For the remainder of this section, we assume \eqref{GeSP}/\eqref{SPRLX} to be augmented with zero Dirichlet boundary conditions for the simplicity of the presentation only as the modification of the numerical method to incorporate periodic boundary conditions is standard. We begin by first presenting the time semi-discrete scheme before moving on to the presentation of the fully-discrete scheme. 

\subsection{Time semi-discrete scheme}\label{TDS}

%Let $H_0^1(\varOmega)$ denote the classical Sobolev space of functions whose derivatives are square integrable on $\varOmega$ and whose traces vanish on the boundary $\partial\varOmega$. 
% Stephen: obvious and uncessary IMO
%We begin by 
We introduce a sequence of $N+1$ \emph{time nodes} $\tau =: t_0 < ... < t_n < \ldots < t_N := T$ %which induce a partition of $(0,T)$ into $N$ open time intervals $I_n := (t_{n-1},t_n)$, $n=1,\ldots,N$. The length of the time interval $I_n$ is called the \emph{time step length} and is given by $k_n := |I_n| = t_n - t_{n-1}$. 
of $[\tau, T]$ and the variable time-steps $k_n:=t_n-t_{n-1}.$
With this notation, our time semi-discrete Besse-style relaxation scheme for \eqref{GeSP} based on \eqref{SPRLX} is defined as follows:  %On the interval $I_n$, 
We seek approximations $(U^{n}, V^{n-\nicefrac{1}{2}},\Phi^{n-\nicefrac{1}{2}}) \in H_0^1(\varOmega)$ to $(u(t_{n}), \V(t_{n-\nicefrac{1}{2}}), \phi(t_{n-\nicefrac{1}{2}})) \in H_0^1(\varOmega)$, $1\le n\le N$, such that
\begin{equation}
\label{CNrelax} \left \{
\begin{aligned}
&\overline{\partial} U^n-\ii p(t_{n-\nicefrac{1}{2}}) \Delta U^{n-\nicefrac{1}{2}}+\ii q(t_{n-\nicefrac{1}{2}})V^{n-\nicefrac{1}{2}}U^{n-\nicefrac{1}{2}}=0,\\
& \Delta V^{n-\nicefrac{1}{2}} =\Phi^{n-\nicefrac{1}{2}}, \\
& \frac{k_{n-1}}{k_{n-1}+k_{n}}\Phi^{n-\nicefrac{1}{2}} = (|U^{n-1}|^{2} -\mu) -  \frac{k_{n}}{k_{n-1}+k_{n}}\Phi^{n-\nicefrac{3}{2}},
\end{aligned}
\right.
\end{equation}
holds,  where we used the notation
\begin{equation}
\label{notation}
t_{n-\nicefrac{1}{2}} := \frac{t_{n-1}+t_{n}}{2}, \quad U^{n - \nicefrac{1}{2}}:=\frac{U^{n-1}+U^{n}}2 ,\quad \overline{\partial} U^n:=\frac{U^{n}-U^{n-1}}{k_{n}}.
\end{equation}
The scheme should be initialized at $n=0;$ a natural, straightforward choice is  $U^0=u_0, \ k_0 = k_1, \ \Phi^{-\nicefrac{1}{2}}=|u_0|^2-\mu.$ In Section \ref{sec:initializ} we discuss some  subtle  computational issues with this initialization, and  derive a modified initialization which addresses them. 
%In the sequel the following easily derived \emph{discrete difference} identity, will be useful
%\begin{equation}
%\label{ddi}
%V^{n-\nicefrac{1}{2}} U^n - V^{n-\nicefrac{3}{2}} U^{n-1} = k_n V^{n-1} \overline{\partial} U^n + \frac{k_n+k_{n-1}}{2} U^{n - \nicefrac{1}{2}} \overline{\partial} V^{n-\nicefrac{1}{2}}
%\end{equation}
%where we have used the notation 
%\begin{equation*}
%V^{n-1}:=\frac{V^{n-\nicefrac{1}{2}} + V^{n-\nicefrac{3}{2}}}{2}, \quad  \overline{\partial} V^{n-\nicefrac{1}{2}}:=\frac{2}{k_{n}+k_{n-1}} (V^{n-\nicefrac{1}{2}} - V^{n-\nicefrac{3}{2}}).
%\end{equation*} 
%If the time step is uniform, $k_n = k, \ n=0,\dots, N$ then the last equation in \eqref{CNrelax} simplifies to 
%\begin{equation*}
%	\Phi^{n-\nicefrac{1}{2}} = 2 |U^{n-1}|^{2} - \Phi^{n-\nicefrac{3}{2}}.
%\end{equation*}
% Stephen: obvious and uncessary IMO
\begin{remark}\upshape
Method \eqref{CNrelax} can be combined with various methods for spatial discretization. Examples include finite differences, in the spirit of \cite{Besse, Zouraris}, spectral methods, or finite elements. In this paper we choose finite elements for the spatial discretization, which allow for spatial adaptivity; this is particularly important for the Schr\"odinger-Poisson system \eqref{GeSP} with applications in cosmology, in which we observe sharply localized features for the density (cf. Section~\ref{sec:cosmoex}).
	
\end{remark}

\subsection{Fully-discrete scheme}

Let $\Th$ be a conforming, shape regular partition of $\varOmega$ consisting of elements $K$ which are either simplices or $d$-dimensional cubes. We then build real/complex finite element spaces over the mesh $\Th$ , denoted by $\Vh(\Th; \mathbb{R})$ and $\Vh(\Th; \mathbb{C})$, respectively, given by
\begin{equation}
\begin{aligned}
\Vh(\Th; \mathbb{R}) & :=\left\{\chi\in C(\bar{\varOmega})\cap H^1_0(\varOmega) :\forall K\in \Th, \  \chi |_K \in \Pb^r(K)\right\}\!, \\
\Vh(\Th; \mathbb{C}) & :=\left\{\chi_R + \ii \chi_I : \chi_R, \chi_I \in \Vh(\Th;\mathbb{R}) \right\}\!,
\end{aligned}
\end{equation}
where $\Pb^r(K)$ denotes the space of polynomials on the element $K$ of total degree $r$ if $K$ is a simplex or of degree $r$ in each variable if $K$ is a $d$-dimensional cube. At each time step $n$, we assume that we have some mesh $\Th^n$ which has been obtained from a previous mesh $\Th^{n-1}$ via a limited number of refinement and/or coarsening operations. We then associate to each time step $n$ the real and complex finite element spaces $\Vh^n(\mathbb{R}) := \Vh(\Th^n;\mathbb{R})$ and $\Vh^n(\mathbb{C}) := \Vh(\Th^n;\mathbb{C})$ over the mesh $\Th^n$.

To characterize the fully-discrete scheme on (possibly) variable finite element spaces, we need to introduce two operators; namely, the $L^2$-projection operator $\Pc_h^n : L^2(\varOmega) \to \Vh^n(\mathbb{C})$ and the discrete Laplacian operator $-\Dh^n : H_0^1(\varOmega) \to \Vh^n(\mathbb{C})$, which are defined implicitly as the solution  of the following variational problems
\begin{align}
& v\mapsto \Pc_h^nv,  &&\hspace{-15mm} \langle \Pc_h^n v, \chi^n \rangle = \langle v , \chi^n \rangle, &\hspace{-15mm}\forall \chi^n \in \Vh^n(\mathbb{R}), \label{L2pr} \\
& v \mapsto -\Dh^n v, &&\hspace{-15mm} \langle -\Dh^n v, \chi^n \rangle = \langle \nabla v , \nabla \chi^n \rangle, &\hspace{-15mm} \forall \chi^n \in \Vh^n(\mathbb{R}), \label{DLapl}
\end{align}
	where $\langle \cdot, \cdot \rangle$ denotes the $L^2$-inner product over $\varOmega$. Note that although the $L^2$-projection/discrete Laplacian may be complex in the above definitions, the test functions always lie in the real finite element space $\Vh^n(\mathbb{R})$. We are now ready to introduce the fully-discrete Besse-style relaxation scheme for \eqref{SP} based on \eqref{SPRLX} which is given as follows: We seek approximations $(U_h^{n}, V_h^{n-\nicefrac{1}{2}}, \Phi_h^{n-\nicefrac{1}{2}})$ $\in \Vh^{n}(\mathbb{C})\!\times\!\Vh^{n}(\mathbb{R})\!\times\!\Vh^{n}(\mathbb{R})$ to $(u(\cdot, t_{n}), \V(\cdot, t_{n-\nicefrac{1}{2}}), \phi(\cdot, t_{n-\nicefrac{1}{2}})) \in H^1_0(\varOmega)$, $1\le n \le N$, such that
\begin{equation}
\label{RCNG} \left \{
\begin{aligned}
&\Pc_h^{n}\bigg[\overline{\partial} U_h^n-\frac{\ii}2 p(t_{n-\nicefrac{1}{2}}) (\Dh^{n-1} U_h^{n-1} +\Dh^{n} U_h^{n})+\ii q(t_{n-\nicefrac{1}{2}}) V_h^{n-\nicefrac{1}{2}}U_h^{n-\nicefrac{1}{2}}\bigg]=0, \\
& \Delta_h^{n} V_h^{n-\nicefrac{1}{2}} = \Phi_h^{n-\nicefrac{1}{2}}, \\
& \frac{k_{n-1}}{k_{n-1}+k_{n}}\Phi_h^{n-\nicefrac{1}{2}} = \Pc_h^{n}\bigg[ (|U_h^{n-1}|^{2} -\mu)- \frac{k_{n}}{k_{n-1}+k_{n}}\Phi_h^{n-\nicefrac{3}{2}}\bigg], 
\end{aligned}
\right.
\end{equation}
where the straightforward initialization would be $U_h^0=\Pc_h^0 u_0$, $k_0 = k_1$ and $\Phi_h^{-\nicefrac{1}{2}}=\Pc_h^0(|u_0|^2-\mu)$. If no mesh change occurs on the time step $n$, i.e., $\Th^{n} = \Th^{n-1}$ then the fully-discrete Besse-style relaxation scheme \eqref{RCNG} can be simplified to
\begin{equation}
\label{RCNGnomeshchange} \left \{
\begin{aligned}
&\overline{\partial} U_h^n-\ii p(t_{n-\nicefrac{1}{2}}) \Dh^{n} U_h^{n-\nicefrac{1}{2}}+\ii q(t_{n-\nicefrac{1}{2}})\Pc_h^{n}(V_h^{n-\nicefrac{1}{2}}U_h^{n-\nicefrac{1}{2}})=0, \\
& \Delta_h^{n} V_h^{n-\nicefrac{1}{2}} = \Phi_h^{n-\nicefrac{1}{2}}, \\
& \frac{k_{n-1}}{k_{n-1}+k_{n}}\Phi_h^{n-\nicefrac{1}{2}} = \Pc_h^{n}(|U_h^{n-1}|^{2}-\mu) - \frac{k_{n}}{k_{n-1}+k_{n}}\Phi_h^{n-\nicefrac{3}{2}}.
\end{aligned}
\right.
\end{equation}

We now look into whether the numerical scheme \eqref{RCNGnomeshchange} satisfies discrete versions of the system's conservation  laws (cf. Lemma \ref{conslaws} for the continuous version). 

The \emph{discrete mass} is the discrete equivalent of the mass  \eqref{cldef1} and is given by
\begin{equation}
\mathcal{M}_h^n := ||U_h^n||^2.
\end{equation}
The discrete mass satisfies an exact equivalent of the conservation of mass law, as is seen in the following:
\begin{proposition}[Discrete Mass Conservation]\label{discmasslemma}
If no mesh change occurs on the time step $n$, i.e. if $\Th^{n} = \Th^{n-1},$ then the solution of the fully-discrete Besse-style relaxation scheme \eqref{RCNGnomeshchange} satisfies
$$\mathcal{M}_h^n = \mathcal{M}_h^{n-1}.$$
Therefore, if no mesh change occurs at all, i.e. if $\Th^{n} = \Th^0,\ \forall 1\le n\le N,$ then 
\begin{equation}
\label{DMC}
\mathcal{M}_h^n  = \mathcal{M}_h^0, \quad 1\le n\le N.
\end{equation}
\end{proposition}
\begin{proof}
We multiply the discrete Schr\"odinger equation in \eqref{RCNGnomeshchange} by $\overline{U}_h^{n-\nicefrac{1}{2}}$, i.e, the complex conjugate of ${U}_h^{n-\nicefrac{1}{2}}$ and integrate over $\varOmega$ to obtain 
\begin{equation}
\int_{\varOmega}\overline{U}_h^{n-\nicefrac{1}{2}}\overline{\partial} U_h^{n} \dif x+\ii p(t_{n-\nicefrac{1}{2}})||\nabla U_h^{n-\nicefrac{1}{2}}||^2+\ii q(t_{n-\nicefrac{1}{2}})\int_{\varOmega} V_h^{n-\nicefrac{1}{2}}|U_h^{n-\nicefrac{1}{2}}|^2 \dif x=0.\end{equation}
The last two terms are purely imaginary so taking real parts and expanding yields
\begin{equation}
\mathcal{M}_h^n - \mathcal{M}_h^{n-1} + \int_{\varOmega}\Rc(U_h^{n}\overline{U}_h^{n-1} - U_h^{n-1}\overline{U}_h^{n}) \dif x = 0.
\end{equation}
The last integral vanishes and so we are left with
\begin{equation}
\mathcal{M}_h^n = \mathcal{M}_h^{n-1},
\end{equation}
as claimed.
\end{proof}
The {\em discrete energy balance} is slightly more sophisticated, as more can be said about how to discretize the energy. In equation \eqref{cldef2}, two expressions for the potential energy were given; both of them come into play, but at the discrete level they are not necessarily identical.
In that context, we will define the {\em discrete kinetic energy} $\E_{k,h}^n,$ and the two {\em discrete versions of the potential energy} $\E_{{\V_1},h}^n, \E_{{\V_2},h}^n$ as follows:
\begin{equation}
\begin{aligned}
\label{DE}
& \E_{k,h}^n  = \|\nabla U_h^n\|^2, \quad  \E_{{\V_1},h}^n= \|\nabla V_h^{n-\nicefrac{1}{2}}\|^2 \quad  \E_{{\V_2},h}^n = -\int_{\varOmega} V_h^{n - \nicefrac{1}{2}}(|U_h^{n}|^2-\mu) \dif x.
\end{aligned}
\end{equation}
What will end up playing the role of {\em the discrete potential energy} would be $\E_{{\V},h}^n:=2\E_{{\V_2},h}^n-\E_{{\V_1},h}^n.$ We are now ready to prove the following
\begin{proposition}[Discrete Energy Balance]\label{discenergylemma}
Assume that the time-step size  remains  constant  between successive time-steps, i.e $k_{n-1}=k_n$, and that  no mesh change occurs on time step $n,$ i.e. $\Th^{n} = \Th^{n-1}.$  Then the solution of the fully-discrete  relaxation scheme \eqref{RCNGnomeshchange} satisfies  the  discrete energy balance law
\begin{align}
\label{localCE}
% \mathcal{E}_h^{n} & = \mathcal{E}_h^{n-1}. 
p(t_{n-\nicefrac{1}{2}})  \overline{\partial}(\E_{k,h}^n) -  \frac{q(t_{n-\nicefrac{1}{2}})}2 \overline{\partial}\left( 2 \E_{{\V_2},h}^n - \E_{{\V_1},h}^n \right)= 0, 
\end{align}
which is the discrete analog of \eqref{cl2}. 

If furthermore $p, q$ are constants,  \eqref{localCE} simplifies to 
\begin{equation}
\label{localBE}
p \mathcal{E}_{k,h}^{n} {\color{black}-}\frac{q}2 (2\mathcal{E}_{\V_2,h}^{n} - \mathcal{E}_{\V_1,h}^{n})  = p \mathcal{E}_{k,h}^{n-1} {\color{black}-}\frac{q}2 (2\mathcal{E}_{\V_2,h}^{n-1} - \mathcal{E}_{\V_1,h}^{n-1}).
\end{equation}
Finally, if in addition the time-step size and spatial mesh do not change over the whole computation, i.e. if $k_n=k$ and $\Th^{n} = \Th$ for all $n=1,2,\dots,N,$ then 
\begin{equation}
\label{DEC}
p \mathcal{E}_{k,h}^{n} +\frac{q}2 (2\mathcal{E}_{\V_2,h}^{n} - \mathcal{E}_{\V_1,h}^{n})  = p \mathcal{E}_{k,h}^{1} +\frac{q}2 (2\mathcal{E}_{\V_2,h}^{1} - \mathcal{E}_{\V_1,h}^{1}).
\end{equation}

%Further, if no mesh change occurs at all, i.e., $\Th^{n} = \ldots = \Th^0$ then 
%\begin{equation}
%\label{ConEne}
%\mathcal{E}_h^n  = \mathcal{E}_h^0.
%\end{equation}
\end{proposition}
  
\begin{proof}
Multiplying the discrete Schr\"odinger equation \eqref{RCNGnomeshchange} by $\overline{\partial}\,\overline{U}_h^n$ and integrating over $\varOmega$ we obtain
\begin{equation*}
\|\overline{\partial} U_h^n\|^2 +\ii p(t_{n-\nicefrac{1}{2}})\!\int_{\varOmega}\nabla U_h^{n - \nicefrac{1}{2}} \bigcdot \nabla \overline{\partial}\,\overline{U}_h^n \dif x + \ii q(t_{n-\nicefrac{1}{2}}) \int_{\varOmega} V_h^{n - \nicefrac{1}{2}}U_h^{n - \nicefrac{1}{2}}\overline{\partial}\,\overline{U}_h^n \dif x = 0.
\end{equation*}
Taking imaginary parts yields
\begin{equation}
\label{ce0}
\frac{p(t_{n-\nicefrac{1}{2}})}{2 k_n}(\|\nabla U_h^{n}\|^2 - \|\nabla U_h^{n-1}\|^2)  + \frac{q(t_{n-\nicefrac{1}{2}})}{2 k_n}\! \int_{\varOmega} V_h^{n - \nicefrac{1}{2}}(|U_h^{n}|^2 - |U_h^{n-1}|^2) \dif x = 0.
\end{equation}
For the second term of \eqref{ce0} we have 
\begin{equation*}
\begin{aligned}
& \int_{\varOmega} V_h^{n - \nicefrac{1}{2}}(|U_h^{n}|^2   - |U_h^{n-1}|^2) \dif x   = \int_{\varOmega} V_h^{n - \nicefrac{1}{2}}\left((|U_h^{n}|^2-\mu) - (|U_h^{n-1}|^2-\mu)\right) \dif x \\ 
  & = \int_{\varOmega} \left(V_h^{n - \nicefrac{1}{2}}(|U_h^{n}|^2-\mu) -  V_h^{n - \nicefrac{3}{2}}(|U_h^{n-1}|^2-\mu)\right) \dif x \\
  & + \int_{\varOmega}\left(V_h^{n - \nicefrac{3}{2}}(|U_h^{n-1}|^2-\mu) - V_h^{n - \nicefrac{1}{2}} (|U_h^{n-1}|^2-\mu)\right) \dif x \\
& = \int_{\varOmega} k_n\ \overline{\partial}\left(V_h^{n - \nicefrac{1}{2}}(|U_h^{n}|^2-\mu) \right) \dif x + \int_{\varOmega} (|U_h^{n-1}|^2-\mu) ( V_h^{n - \nicefrac{3}{2}} - V_h^{n - \nicefrac{1}{2}}) \dif x
\end{aligned}
\end{equation*}
Hence \eqref{ce0} becomes 
\begin{equation}
\label{ce1}
\begin{aligned}
\frac{p(t_{n-\nicefrac{1}{2}})}{2} &\overline{\partial}(\|\nabla U_h^{n}\|^2) + \\
+& \frac{q(t_{n-\nicefrac{1}{2}})}{2 k_n}\left(\!\int_{\varOmega}\!\! k_n \overline{\partial}\left(V_h^{n - \nicefrac{1}{2}}(|U_h^{n}|^2-\mu) \right) \dif x \!+ \!\! \int_{\varOmega}\!( |U_h^{n-1}|^2-\mu) ( V_h^{n - \nicefrac{3}{2}} - V_h^{n - \nicefrac{1}{2}}) \dif x\!\!\right)\!\!=\!0.
\end{aligned}
\end{equation} 
Using \eqref{RCNG}(c), \eqref{RCNG}(b), and finally integration by parts,  for the third term of \eqref{ce1}  we obtain,
\begin{equation}
\label{ce2}
\begin{aligned}
&\int_{\varOmega}\!( |U_h^{n-1}|^2-\mu) ( V_h^{n - \nicefrac{3}{2}} - V_h^{n - \nicefrac{1}{2}}) \dif x=\int_{\varOmega}\!\Pc_h^{n}\Big( |U_h^{n-1}|^2-\mu\Big) ( V_h^{n - \nicefrac{3}{2}} - V_h^{n - \nicefrac{1}{2}}) \dif x\\
& =
\int_{\varOmega} \left( \frac{k_{n-1}}{k_{n-1}+k_{n}}\Phi_h^{n-\nicefrac{1}{2}} + \frac{k_{n}}{k_{n-1}+k_{n}}\Phi_h^{n-\nicefrac{3}{2}}\right) ( V_h^{n - \nicefrac{3}{2}} - V_h^{n - \nicefrac{1}{2}}) \dif x \\ 
& =  \int_{\varOmega} \left( \frac{k_{n-1}}{k_{n-1}+k_{n}}\Delta_h^{n} V_h^{n-\nicefrac{1}{2}}+ \frac{k_{n}}{k_{n-1}+k_{n}}\Delta_h^{n} V_h^{n-\nicefrac{3}{2}}\right) ( V_h^{n - \nicefrac{3}{2}} - V_h^{n - \nicefrac{1}{2}}) \dif x  \\
& = \int_{\varOmega} \left( \frac{k_{n-1}}{k_{n-1}+k_{n}}\nabla V_h^{n-\nicefrac{1}{2}}+ \frac{k_{n}}{k_{n-1}+k_{n}}\nabla V_h^{n-\nicefrac{3}{2}}\right) ( \nabla V_h^{n - \nicefrac{1}{2}} - \nabla V_h^{n - \nicefrac{3}{2}}) \dif x.%\\
%& = \frac12\int_{\varOmega} \left( |\nabla V_h^{n-\nicefrac{1}{2}}|^2 - |\nabla V_h^{n-\nicefrac{3}{2}}|^2\right) \dif x  = \frac12 k_n \overline{\partial}(\|\nabla V_h^{n-\nicefrac{1}{2}}\|^2),
\end{aligned}
\end{equation}
Since the time-step size is constant between successive time-steps we readily obtain that
\begin{equation*}
%\begin{aligned}
\int_{\varOmega}\!( |U_h^{n-1}|^2-\mu) ( V_h^{n - \nicefrac{3}{2}} - V_h^{n - \nicefrac{1}{2}}) \dif x= \frac12\int_{\varOmega} \left( |\nabla V_h^{n-\nicefrac{1}{2}}|^2 - |\nabla V_h^{n-\nicefrac{3}{2}}|^2\right) \dif x  = \frac12 k_n \overline{\partial}(\|\nabla V_h^{n-\nicefrac{1}{2}}\|^2)
%\end{aligned}
\end{equation*}
%where the last equality holds true since the time-step is constant   between successive time-steps, $k_n=k_{n-1}.$  
Using the above to reformulate the last term of \eqref{ce1} yields
\begin{equation*}
\frac{p(t_{n-\nicefrac{1}{2}})}{2} \overline{\partial}(\|\nabla U_h^{n}\|^2) + \frac{q(t_{n-\nicefrac{1}{2}})}{2k_n} \left(k_n\int_{\varOmega}  \overline{\partial}\left(V_h^{n - \nicefrac{1}{2}}(|U_h^{n}|^2-\mu) \right) \dif x + \frac{k_n}{2}\overline{\partial}(\|\nabla V_h^{n-\nicefrac{1}{2}}\|^2)  \right)= 0 .
\end{equation*}
The result then follows from the definition of the discrete energies $\E_{k,h}^n, \ \E_{{\V_1},h}^n, \ \E_{{\V_2},h}^n$, \eqref{DE}.
\end{proof}

\begin{remark}[Discrete balance laws \& periodic boundary conditions]\upshape
Propositions~\ref{discmasslemma} \& \ref{discenergylemma}	remain valid for the fully discrete scheme corresponding to \eqref{GeSP}/\eqref{SPRLX} with periodic boundary conditions, as long as the finite element space $\Vh(\Th; \mathbb{C})$ is appropriately equipped with periodic boundary conditions instead of homogeneous Dirichlet.
\end{remark}

An important question that raises here is what happens to the discrete energy balance in the case of variable time-steps. In particular, by how much does it fail to satisfy \eqref{localCE}? We answer this in the next proposition:
%{\color{black} What happens to the energy balance when you change the time-step size?

\begin{proposition}	[Discrete Energy Balance \& Variable time-steps]\label{lm:jump}
Assume that  no mesh change occurs on time step $n,$ i.e. $\Th^{n} = \Th^{n-1}.$  Then the solution of the fully-discrete  relaxation scheme \eqref{RCNGnomeshchange} satisfies  the following %discrete energy balance 
\begin{align}
\label{localCE2}
% \mathcal{E}_h^{n} & = \mathcal{E}_h^{n-1}. 
p(t_{n-\nicefrac{1}{2}})  \overline{\partial}(\E_{k,h}^n) -  \frac{q(t_{n-\nicefrac{1}{2}})}2 \overline{\partial}\left( 2 \E_{{\V_2},h}^n - \E_{{\V_1},h}^n \right)+\frac{q(t_{n-1/2})}{\color{black}4}\frac{k_{n-1}-k_n}{k_{n-1}+k_n}k_n\|\bar\p\nabla V_h^{n-1/2}\|^2= 0.
\end{align}
\end{proposition}

\begin{proof}
For this case, \eqref{ce1} and \eqref{ce2} are still valid. From \eqref{ce2} we obtain
\begin{equation*}
%\label{ce2}
\begin{aligned}
&\int_{\varOmega}\!( |U_h^{n-1}|^2-\mu) ( V_h^{n - \nicefrac{3}{2}} - V_h^{n - \nicefrac{1}{2}}) \dif x=\int_{\varOmega}\!\Pc_h^{n}\Big( |U_h^{n-1}|^2-\mu\Big) ( V_h^{n - \nicefrac{3}{2}} - V_h^{n - \nicefrac{1}{2}}) \dif x\\=
&\int_{\varOmega} \left( \frac{k_{n-1}}{k_{n-1}+k_{n}}\Phi_h^{n-\nicefrac{1}{2}} + \frac{k_{n}}{k_{n-1}+k_{n}}\Phi_h^{n-\nicefrac{3}{2}}\right) ( V_h^{n - \nicefrac{3}{2}} - V_h^{n - \nicefrac{1}{2}}) \dif x\\
=&\frac12 k_n \overline{\partial}(\|\nabla V_h^{n-\nicefrac{1}{2}}\|^2)+\frac{k_{n-1}-k_n}{2(k_{n-1}+k_n)}\int_{\varOmega}(\nabla V_h^{n-1/2}-\nabla V_h^{n-3/2})^2\dif x
\\
=&\frac12 k_n \overline{\partial}(\|\nabla V_h^{n-\nicefrac{1}{2}}\|^2)+\frac{k_{n-1}-k_n}{2(k_{n-1}+k_n)}k_n^2\|\bar\p\nabla V_h^{n-1/2}\|^2.
\end{aligned}
\end{equation*}	
Using the above to reformulate the last term of \eqref{ce1} we obtain \eqref{localCE2}.
\end{proof}

\begin{remark}\upshape
Clearly, if $k_{n-1}=k_n$ \eqref{localCE2} reduces to \eqref{localCE}. For variable time-steps the remainder term 	$\dfrac{q(t_{n-1/2})}4\dfrac{k_{n-1}-k_n}{k_{n-1}+k_n}k_n\|\bar\p\nabla V_h^{n-1/2}\|^2$ is expected to be of first order in time; in other words, the energy balance law \eqref{localCE} fails to be satisfied by order $\mathcal{O}(k)$ (where $k=\max_{1\le n\le N-1}k_n$), every time we change the time-step. %The result in Proposition~\ref{lm:jump} is further investigated numerically in Section \ref{sec:numjump}.
\end{remark}

\section{Implementation} \label{sec:implemntation}

In this section, we discuss the practicalities of implementing the Besse-style relaxation scheme \eqref{RCNG} for the numerical solution of the Schr\"odinger-Poisson system \eqref{SP}.

\subsection{Initialization} \label{sec:initializ}
In Section \ref{timediscr}, the straightforward initialization $U_h^0=\Pc_h^0 u_0,$ $k_0 = k_1$ and $\Phi_h^{-\nicefrac{1}{2}}=\Pc_h^0(|u_0|^2-\mu)$ was presented. This is a simple, viable choice, and we observe numerically that the obtained numerical solution $U^n$ is a second order approximation in time to $u(t_n)$. However, in the same computations we observe that $\Phi^{n-\frac{1}2}$ is only a first order approximation in time to $\V(t_{n-1/2})$. This is something also observed in \cite{Zouraris}. Thus we look for a modified initialization under which both $U^n$ and $\Phi^{n-\frac{1}2}$ will be seen numerically to be of second order in time.

One way to do this is to define
 $\Phi^{\nicefrac{1}{2}}_{h,\text{old}}$ according to the straightforward initialization, i.e. $\Phi^{\nicefrac{1}{2}}_{h,\text{old}}=\Pc_h^0(|u_0|^2-\mu)$. $\Phi^{\nicefrac{1}{2}}_{h, \text{old}}$ is then used in the numerical scheme \eqref{RCNG} to calculate initial approximations for the potential and wavefunction which we denote by $\widetilde{V}^{\nicefrac{1}{2}}_h$ and $\widetilde{U}^1_h$, respectively. The initial approximation to the wavefunction, $\widetilde{U}_h^1$, is then used to update the estimate for the nonlinearity on the first time step. In particular, the coefficients $\widehat{\Phi}_h^{\nicefrac{1}{2}}$ are chosen to satisfy
\begin{equation}\label{eq:phi0_1}M\widehat{\Phi}_h^{\nicefrac{1}{2}} = \bigg(\bigg\langle\frac{1}{2}(|\widetilde{U}_h^1|^2-\mu) + \frac{1}{2}{\Phi}_{h,\text{old}}^{\nicefrac{1}{2}},\varphi^1_j \bigg \rangle_{\!\!j=1,\ldots,\mathcal{N}_{1}} \bigg).
\end{equation}
This is equivalent to choosing a modified 
\begin{equation}\label{eq:phi0_}
\Phi^{-\nicefrac{1}{2}}_h=\frac{3}2\left(|U^0_h|^2-\mu\right) - \frac{1}2\left( |\widetilde{U}^1_h|^2-\mu\right)
\end{equation}
to be used in \eqref{RCNG}(c). 
%where 
%\begin{equation}
%\label{RCNG_0} \left \{
%\begin{aligned}
%&\Pc_h^{0}\bigg[\frac{\widetilde{U}_h^1 - {U}_h^0}{k_0}-\frac{\ii}2 p(t_{-\nicefrac{1}{2}}) (\Dh^{-1} U_h^{0} +\Dh^{0} \widetilde{U}_h^{1})+\ii q(t_{-\nicefrac{1}{2}}) V_h^{-\nicefrac{1}{2}}( U_h^{0} + \widetilde{U}_h^{1})\bigg]=0, \\
%& \Delta_h^{0} V_h^{-\nicefrac{1}{2}} =  \Pc_h^{n}\bigg[ (|U_h^{0}|^{2} -\mu)\bigg]. 
%\end{aligned}
%\right.
%\end{equation}
Using this initialization we observe numerically second order in time for both $U^n,$ $\Phi^{n-\frac{1}2}.$ An analogous initialization can be found in \cite{Zouraris} for a Besse-type relaxation finite difference scheme and the semilinear parabolic equation.

\subsection{Solving for the nonlinearity}

Solving for the nonlinearity is a standard finite element problem, i.e., we are seeking a vector of coefficients $\widehat{\Phi}_h^{n-\nicefrac{1}{2}}$ for ${\Phi}_h^{n-\nicefrac{1}{2}}$ such that
$${\Phi}_h^{n-\nicefrac{1}{2}} = \sum_{j = 1}^{\mathcal{N}_{n}} \widehat{\Phi}_{h,j}^{n-\nicefrac{1}{2}}\varphi^{n}_j,$$
where $\varphi^{n}_j$, $j = 1, \ldots, \mathcal{N}_{n} = \text{dim}(\Vh^{n}(\mathbb{R}))$ are (real) finite element basis functions. From \eqref{RCNG}, we see that the vector of coefficients $\widehat{\Phi}_h^{n-\nicefrac{1}{2}}$ must satisfy
$$\frac{k_{n-1}}{k_{n-1}+k_{n}}M\widehat{\Phi}_h^{n-\nicefrac{1}{2}} = \bigg(\bigg\langle|U_h^{n-1}|^2 - \frac{k_{n}}{k_{n-1}+k_{n}}{\Phi}_h^{n-\nicefrac{3}{2}},\varphi_j^{n} \bigg \rangle_{\!\!j=1,\ldots,\mathcal{N}_{n}} \bigg),$$
where $M$ is the \emph{mass matrix} given by
$$M_{ij} = \int_{\varOmega} \varphi^{n}_i \varphi^{n}_j \, \dif x.$$

\subsection{Solving for the potential}

Solving for the potential is also fairly routine; if we introduce the vector of coefficients $\widehat{V}_h^{n-\nicefrac{1}{2}}$ such that
$$V_h^{n-\nicefrac{1}{2}} = \sum_{j = 1}^{\mathcal{N}_{n}} \widehat{V}_{h,j}^{n-\nicefrac{1}{2}}\varphi^{n}_j,$$
then \eqref{RCNG} implies that $\widehat{V}_h^{n-\nicefrac{1}{2}}$ must satisfy 
$$L\widehat{V}_h^{n-\nicefrac{1}{2}} = \bigg(\langle \Phi_h^{n-\nicefrac{1}{2}}, \varphi^{n}_j \rangle_{j=1,\ldots,\mathcal{N}_{n}}\bigg),$$
where $L$ is the \emph{stiffness matrix} given by 
$$L_{ij} = \int_{\varOmega} \nabla\varphi^{n}_i \bigcdot \nabla\varphi^{n}_j \, \dif x.$$
For $V_h^0$, the coefficients $\widehat{V}_h^{0}$ are chosen to satisfy the matrix-vector system
\begin{equation}\label{eq:v00}
L\widehat{V}_h^{0} = \bigg(\langle \Phi^{\nicefrac{1}{2}}_{h,\text{old}}, \varphi^0_j \rangle_{j=1,\ldots,\mathcal{N}_0}\bigg).
\end{equation}
Another question now is how to use the nodal values $V_h^{n-1/2}$ in order to obtain optimal (second) order approximations $V_h^n$ to $\V(t_n)$ (note that that the obtained by the method approximations $V_h^{n-1/2}$ are approximations to $\V(t_{n-1/2})$, i.e., are approximations at the middle nodal points $t_{n-1/2}$ and not at the nodes $t_n$).
%An important question now arises: how do we extend these nodal values so that $V_h(t)$ is of optimal order in time? On the first interval, 
To obtain $V_h^1$, the most obvious choice is to linearly extrapolate from $V_h^0$ and $V_h^{\nicefrac{1}{2}}$ to obtain $V_h^1$, i.e., we set $V_h^1 = 2V_h^{\nicefrac{1}{2}} - V_h^0$. %We can then define
%$$V_h(t) := \bigg(\frac{t-t_0}{k_1}\bigg)V_h^1 + \bigg(\frac{t_1-t}{k_1}\bigg)V_h^0, \qquad t \in \bar{I}_1.$$
It is tempting to continue to iterate this procedure in order to compute the remaining nodal values, however,  for $n\neq 0$, %the first interval, 
the value $V_h^{n}$ is an \emph{extrapolated quantity} (in contrast to $V^0_h$ which is computed according to \eqref{eq:v00}). Attempting to calculate $V_h^{n}$ by extrapolating through $V_h^{n-1}$ (an extrapolated point) and $V_h^{n-\nicefrac{1}{2}}$ (a computed point) is therefore an unstable procedure which oscillates out of control. To avoid this, we instead calculate $V_h^{n}$ by linearly extrapolating from the computed values $V_h^{n-\nicefrac{1}{2}}$ and $V_h^{n-\nicefrac{3}{2}}$; a simple calculation yields
$$V_h^{n} =  V_h^{n-\nicefrac{1}{2}} + \frac{k_{n}}{k_{n-1}+k_{n}}(V_h^{n-\nicefrac{1}{2}} - V_h^{n-\nicefrac{3}{2}}).$$
%A variety of options now exist as to how we can define $V_h(t)$ on $\bar{I}_{n}$ using the three nodal values $V_h^{n-1}$, $V_h^{n-\nicefrac{1}{2}}$ and $V_h^{n}$. We opt for the ``standard'' approach of defining $V_h(t)$ to be the linear interpolant through the extremal nodes, i.e., we set
%$$V_h(t) := \bigg(\frac{t_{n}-t}{k_{n}}\bigg)V_h^{n-1}+\bigg(\frac{t-t_{n-1}}{k_{n}}\bigg)V_h^{n}, \qquad t \in \bar{I}_{n}.$$
%We remark, however, that this approach does have a drawback: in general, $V_h(t_{n-\nicefrac{1}{2}}) \neq V_h^{n-\nicefrac{1}{2}}$ (although this does not affect the order of $V_h$).

\subsection{Solving for the wavefunction}

In the case of the wavefunction, we are seeking a vector of real coefficients $\widehat{U}^{n}_{h,R}$ and a vector of imaginary coefficients $\widehat{U}^{n}_{h,I}$ such that
$$U_h^{n} = U^{n}_{h,R} + iU^{n}_{h,I} =  \sum_{j=1}^{\mathcal{N}_{n}} \big(\widehat{U}^{n}_{h,R,j} + i\widehat{U}^{n}_{h,I,j}\big)\varphi^{n}_j.$$
Here, as before, $\varphi^{n}_j$, $j = 1, \ldots, \mathcal{N}_{n} = \text{dim}(\Vh^{n}(\mathbb{R}))$ are \emph{real} finite element basis functions which form a basis for $\Vh^{n}(\mathbb{R})$. Then \eqref{RCNG} implies that the coefficient vectors must satisfy the block matrix-vector system
\begin{equation}
\begin{aligned}
\label{wavefunctionsystem}\displaystyle&\begin{pmatrix}M &  \frac{p(t_{n-\nicefrac{1}2}) k_{n}}{2}L -\frac{q(t_{n-\nicefrac{1}2}) k_{n}}{2}P  \\ -\frac{p(t_{n-\nicefrac{1}2}) k_{n}}{2}L +\frac{q(t_{n-\nicefrac{1}2}) k_{n}}{2}P & M \end{pmatrix}\!\!\begin{pmatrix}\widehat{U}^{n}_{h,R} \\ \widehat{U}^{n}_{h,I} \end{pmatrix} =\\ & \begin{pmatrix}\big\langle U_{h,R}^{n-1} - \frac{p(t_{n-\nicefrac{1}2}) k_{n}}{2}\Delta_h^{n-1}U^{n-1}_{h,I} - \frac{q(t_{n-\nicefrac{1}2}) k_{n}}{2}V_h^{n-\nicefrac{1}{2}}U^{n-1}_{h,I}, \varphi_j^{n} \big\rangle_{j=1,\ldots,\mathcal{N}_{n}} \\ \big\langle U_{h,I}^{n-1} + \frac{p(t_{n-\nicefrac{1}2}) k_{n}}{2}\Delta_h^{n-1}U^{n-1}_{h,R} + \frac{q(t_{n-\nicefrac{1}2}) k_{n}}{2}V_h^{n-\nicefrac{1}{2}}U^{n-1}_{h,R}, \varphi_j^{n} \big\rangle_{j=1,\ldots,\mathcal{N}_{n}} \end{pmatrix},
\end{aligned}
\end{equation}
where $\Delta_h^{n-1}$ is the discrete Laplacian operator (see \eqref{DLapl}), $L$ is the stiffness matrix  and $P$, the matrix associated with the potential term, is given by
$$P_{ij} = \int_{\varOmega} V_h^{n-\nicefrac{1}{2}}\varphi^{n}_i\varphi^{n}_j \,\dif x.$$
As is standard, one can extend the nodal values of the wavefunction to a function $U_h(t)$ on the whole interval via linear interpolation, viz.,
$$U_h(t) := \bigg(\frac{t_{n}-t}{k_{n}}\bigg)U_h^{n-1}+\bigg(\frac{t-t_{n-1}}{k_{n}}\bigg)U_h^{n}, \qquad t \in [t_{n-1},t_n].$$
As a side, we note that \eqref{wavefunctionsystem} can be solved far more efficiently if no mesh change has occurred Indeed, in this case the system \eqref{RCNG} can be rewritten to solve for the half point $U_h^{n-\nicefrac{1}{2}}$ \eqref{RCNGnomeshchange} resulting in the block matrix-vector system
$$\displaystyle\begin{pmatrix}M & - \frac{p(t_{n-\nicefrac{1}2}) k_{n}}{2}L -\frac{q(t_{n-\nicefrac{1}2}) k_{n}}{2}P  \\ \frac{p(t_{n-\nicefrac{1}2}) k_{n}}{2}L +\frac{q(t_{n-\nicefrac{1}2}) k_{n}}{2}P & M \end{pmatrix}\!\!\begin{pmatrix}\widehat{U}^{n-\nicefrac{1}{2}}_{h,R} \\ \widehat{U}^{n-\nicefrac{1}{2}}_{h,I} \end{pmatrix} = \begin{pmatrix} M & 0 \\ 0 & M\end{pmatrix}\!\!\begin{pmatrix} \widehat{U}^{n-1}_{h,R} \\ \widehat{U}^{n-1}_{h,I}\end{pmatrix}.$$
The nodal value coefficients $\widehat{U}_h^{n}$ can then be recovered via $\widehat{U}_h^{n} = 2\widehat{U}_h^{n-\nicefrac{1}{2}} - \widehat{U}_h^{n-1}$. 

\section{Numerical Experiments}
We perform four sets of numerical experiments. First of all we  apply the new numerical method \eqref{RCNG} to some relatively simple problems in order to confirm numerically the rate of convergence. In addition, we verify the validity of the discrete conservation laws, for both time-independent and time-dependent coefficients $p$ and $q$. It must be noted that in what follows we use the modified initialization discussed in Section \ref{sec:initializ}. \RED{We also study how variable time-step affects the conservation of mass and balance of energy.} 

Moreover, we apply \eqref{RCNG} to an example with time-dependent coefficients, periodic boundary conditions and singular features (``sine wave collapse'') which arises in cosmology.  In that context, the semiclassical Schr\"odinger-Poisson system \eqref{SP} is used as a lower dimensional analogue of the  Vlasov-Poisson system \cite{KVS}. The numerical results reported in this section take place in two spatial dimensions and utilize a C++ code based on the {\tt deal.II} finite element library \cite{BHK}.

\subsection{Experimental order of convergence}\label{sec:eoc}
To verify the experimental order of convergence of the numerical method, we apply the classical method of manufactured solutions, i.e., we choose a wavefunction $u(\xv,t):\varOmega\times[0, T]\to \Cb $ and a potential $\V(\xv,t):\varOmega\times[0, T]\to \Rb$ such that  \eqref{SP} is satisfied (with the inclusion of appropriate right-hand sides). Note that in this case $\mu=0$ and homogeneous Dirichlet boundary conditions are used. Moreover we set $\varOmega = (-1,1)^2 \subset \Rb^2$ and consider uniform partitions $\Th$ of $\varOmega$ consisting of squares with sides of length $h$. For simplicity, we set the PDE coefficients to $\alpha=\beta=\ep=1${\color{black}; this leads to $p(t)=\frac12 , \ q(t)=1, \forall t $. The initial time is set to be $\tau=0$ and the final time is given by $T = 1$. The time interval $(\tau,T)$ is} subdivided into uniform intervals of time step length $k = T/N > 0$. % \RED{With choices,  we have $p(t)=\frac12 , \ q(t)=1, \forall t $.} 
 Then, we choose the right-hand sides such that the exact solution to \eqref{SP} is given by 
\begin{equation}\label{exactSP}
\V(x,y,t) = e^{-t}\sin(\pi(x^2-1)(y^2-1)), \qquad\qquad\qquad u(x,y,t) = (1+i)\V(x,y,t).
\end{equation}
The errors are then measured in the $\displaystyle L^{\infty}(L^2)$ norm and we expect that 
$$
e(u;h,k) := \max_{0\le n\le N} \|u(\cdot, t_n) - U_h^n\| = \mathcal{O}(k^2 + h^{r+1}), \quad e({\V};h,k) := \max_{0\le n\le N} \|\V(\cdot, t_n) - V_h^n \| = \mathcal{O}(k^2 + h^{r+1}), 
$$ 
where $r$ is the polynomial degree of the spatial finite element space $\Vh(\Th;\mathbb{R})$.

To compute the spatial convergence rate, we take a large number of time steps, $N=2000$ i.e. $k = 5\cdot 10^{-4}$, so that the temporal part of the error is negligible. We then compute the spatial experimental order of convergence by performing two different realizations with the mesh sizes $h_1$ and $h_2$ and computing
$$
\text{Rate} := \frac{\log(e(\cdot\,;h_1,k)) - \log(e(\cdot\,;h_2,k))}{\log(h_1) - \log(h_2)}.
$$
In Table \ref{eoch}, the spatial experimental orders of convergence are displayed for $r = 1$ and $r = 2$. The optimal rate of convergence is observed in both cases (two for $r=1$ and three for $r=2$) thus validating the claimed spatial accuracy of our new numerical method \eqref{RCNG}.
\begin{table}[htp]
\caption{Spatial experimental orders of convergence for $u, \V$.}\vspace{-2mm}
\begin{center}
\begin{tabular}{|c|cc|cc||cc|cc||} \hline
&  \multicolumn{4}{c||}{$r=1$} & \multicolumn{4}{c||}{$r=2$}  \\ \hline 
$h$ & $e(u;k,h)$ & Rate  &  $e(\V;k,h)$ & Rate & $e(u;k,h)$ & Rate  &  $e(\V;k,h)$ & Rate\\ \hline 
0.250000	   &  2.60203e-1   &	     -      & 	1.36736e-1 & -	       &  1.54310e-2	&      -       &	 8.50485e-3    & -\\ \hline
0.125000	   &  6.58945e-2   &	1.981  & 	3.29791e-2 & 2.052 &  2.19359e-3	& 2.814	   &    1.31987e-3	& 2.688 \\ \hline
0.062500	   &  1.68103e-2   &   1.971  & 	8.23356e-3 & 2.002 &  2.54266e-4	& 3.109	   &    1.71600e-4	& 2.943\\ \hline
0.031250	& 4.22146e-3   & 	1.994  &	2.05895e-3 & 2.000 &  3.12572e-5	& 3.024	   &    2.16460e-5	& 2.987 \\ \hline
0.015625	& 1.05487e-3   & 	2.001  &	5.14783e-4 & 2.000 &  3.85637e-6	& 3.019	   &    2.71179e-6	& 2.997 \\ \hline
%0.250000 &	2.35047e-2	  &            &		1.91319e-2     & 	         \\ \hline
%0.125000 & 	5.96086e-3	  & 1.979  &		4.74129e-3	 &  2.013   \\ \hline
%0.062500 &  1.50954e-3	  & 1.981	&   	1.18408e-3	 &  2.002   \\ \hline
%0.031250 &	3.80764e-4	  & 1.987	&	    2.96019e-4	 &  2.000   \\ \hline
%0.015625 &   9.53201e-5	  & 1.998	&	    7.40047e-5	 &  2.000   \\ \hline
\end{tabular}
\end{center}
\label{eoch}
\end{table}%

For the temporal error rate, we take a large polynomial degree, $r = 3$, in order to minimize the spatial error over the uniform spatial mesh $\mathcal{T}_h$ of mesh size $h = 0.0625$. We then compute the temporal experimental order of convergence by performing two different realizations with the time step lengths $k_1$ and $k_2$ and computing
$$
\text{Rate} := \frac{\log(e(\cdot\,;h,k_1)) - \log(e(\cdot\,;h,k_2))}{\log(k_1) - \log(k_2)}.
$$
The results, given in Table \ref{eock}, confirm that our proposed numerical method \eqref{RCNG}  is of order two in time for both the wavefunction $u$ and the potential $\V$.
\begin{table}[htp]
\caption{Temporal experimental orders of convergence for $u, \V$.}\vspace{-2mm}
\begin{center}
\begin{tabular}{|c|cc|cc|} \hline
$k$ & $e(u;k,h)$ & Rate  &  $e(\V;k,h)$ & Rate\\ \hline 
%0.10	   &  2.13604e-3 &	     -      & 	5.34525e-3 & -	       \\ \hline
0.04	   &  3.72233e-4  &	     -    & 	9.60801e-4 & -	       \\ \hline
0.02	   &  9.49430e-5   &	1.971  & 	2.51017e-4 & 1.936 \\ \hline
0.01	   &  2.39046e-5   &   1.990  & 	6.41950e-5 & 1.967 \\ \hline
\end{tabular}
\end{center}
\label{eock}
\end{table}%

\subsection{Discrete Conservation laws} \label{sec:conslnum}
In this example, we investigate the behaviour of our numerical scheme \eqref{RCNG} for system \eqref{SP}, with respect to the mass conservation  \eqref{cl1} and energy balance \eqref{cl2} in two different cases: a) with constant coefficients $p(t)$ and $q(t)$ and b) with variable coefficients $p(t)$ and $q(t)$. 
We take $\varOmega = (-1,1)^2$ which we discretize with linear finite elements over a uniform grid $\Th$ consisting of squares with sides of length $h=0.015625$. We take  $\tau=0$ and  final time  $T = 3$.  We use $N = 3000$ time steps giving a time step length of $k=10^{-3}$ and we choose 
$\alpha = \beta = 5$.   The initial condition is taken to be
\begin{equation}
\label{CLIC}
u(x,y,0) = \left(\sin(\frac{x}{\pi})+i\cos(\frac{y}{\pi}) \right)(1-x^2)(1-y^2),
\end{equation}
which vanishes along the boundary of $\Omega$ ($\mu=0$ in \eqref{GeSP}). 

\subsubsection{Constant Coefficients $p(t), \ q(t)$} In this case we take $p(t)=\frac{\ep}{50}, \ q(t)=\frac{1}{\ep}, \forall t$ and we expect that both mass and energy are conserved at the discrete level. We allow $\ep$  to vary in order to analyze how this affects the errors in the conservation laws.  {\color{black}We then compute the global conservation law errors \eqref{DMC},  \eqref{DEC}, given by
\begin{align}
\mathcal{M}_e^n & :=\left|\mathcal{M}_h^n  - \mathcal{M}_h^0\right|,  \label{dMC}\\
 \E_{e,gl}^n & := \bigg|\left( p \mathcal{E}_{k,h}^{n} +\frac{q}2 (2\mathcal{E}_{\V_2,h}^{n} - \mathcal{E}_{\V_1,h}^{n})\right)  -  \left(p \mathcal{E}_{k,h}^{1} +\frac{q}2 (2\mathcal{E}_{\V_2,h}^{1} - \mathcal{E}_{\V_1,h}^{1})\right)  \bigg |. \label{dEC}
%:= \left|\ \! \|U_h^n\|^2 - \|u_0\|^2\ \! \right| , \label{dcl1} \\
% \E_e^n & :=\bigg|\left(\frac{\ep^2}{\alpha^2}||\nabla U_h^n||^2 + \beta\|\nabla V_h^{n-\nicefrac{1}{2}}\|^2 + 2\beta\!\int_{\varOmega} V_h^{n-\nicefrac{1}{2}}|U_h^n|^2 \dif x \right) \\
%\notag& \qquad - \left(\frac{\ep^2}{\alpha^2}||\nabla U^1_h||^2 + \beta||\nabla V^{\nicefrac{1}2}_h||^2 + 2\beta\!\int_{\varOmega} V^{\nicefrac{1}2}_h|U^1_h|^2 \dif x \right) \bigg |,
\end{align}
(Note that this definition of global energy error only applies to constant coefficients.)
}
From Table \ref{invs}, we observe that the density and the energy are conserved to double precision accuracy for all values of $\varepsilon$ as expected. {\em Note that each row is roughly 1000 time-steps after the previous one.}
\begin{table}[htp]
\caption{Errors in the conservation laws : $p(t),$ $q(t)$ constant }\vspace{-5mm}
\begin{center}
\begin{tabular}{|c|c|c||c|c||c|c||} \hline
   &  \multicolumn{2}{c||}{$\varepsilon=1$} & \multicolumn{2}{c||}{$\varepsilon=0.1$}  & \multicolumn{2}{c||}{$\varepsilon=0.01$}  \\ \hline
$t_n$  & 	$ \mathcal{M}_e^n$    & 	$ \E_{e,gl}^n$   & 	$ \mathcal{M}_e^n$     & 	$ \E_{e,gl}^n$  &  $ \mathcal{M}_e^n$    &  $ \E_{e,gl}^n$ \\ \hline 
0	&	4.55e-15	&	2.39e-16	&		7.22e-16	&	2.58e-15	&		3.94e-15 	&	9.57e-15	\\ \hline 
1	&	2.06e-14	&	3.29e-16	& 		4.11e-15	&	1.39e-15	&		3.55e-15	&	1.60e-14	\\ \hline
2	&	4.33e-14	&	1.75e-16	&		1.66e-15	&	2.36e-15	&		8.55e-15	&	1.54e-14	\\ \hline
3	&	5.97e-14	&	3.03e-16	&		7.32e-15	&	2.01e-15	&		1.44e-14	&	2.56e-14	\\ \hline
\end{tabular}
\end{center}
\label{invs}%\vspace{-10mm}
\end{table}

\subsubsection{Variable Coefficients $p(t),$ $q(t)$} For the variable coefficient case we take $p(t) = \frac{\ep}{50}t, \ q(t) =\frac{1}{\ep} t^{\frac12}, \forall t$. {\color{black}Here we have to introduce the local error for the discrete energy balance law, namely
\begin{equation}
\begin{aligned}
 \E_{e,loc}^n  :=    \bigg|&\left( p(t_{n-\nicefrac{1}{2}})  \E_{k,h}^n -  \frac{q(t_{n-\nicefrac{1}{2}})}2 \left( 2 \E_{{\V_2},h}^n - \E_{{\V_1},h}^n \right) \right) - 
\\ &\hspace{2.25cm}\left( p(t_{n-\nicefrac{3}{2}})  \E_{k,h}^{n-1} -  \frac{q(t_{n-\nicefrac{3}{2}})}2 \left( 2 \E_{{\V_2},h}^{n-1} - \E_{{\V_1},h}^{n-1} \right) \right)  \bigg |. \label{dECloc}
%:= \left|\ \! \|U_h^n\|^2 - \|u_0\|^2\ \! \right| , \label{dcl1} \\
% \E_e^n & :=\bigg|\left(\frac{\ep^2}{\alpha^2}||\nabla U_h^n||^2 + \beta\|\nabla V_h^{n-\nicefrac{1}{2}}\|^2 + 2\beta\!\int_{\varOmega} V_h^{n-\nicefrac{1}{2}}|U_h^n|^2 \dif x \right) \\
%\notag& \qquad - \left(\frac{\ep^2}{\alpha^2}||\nabla U^1_h||^2 + \beta||\nabla V^{\nicefrac{1}2}_h||^2 + 2\beta\!\int_{\varOmega} V^{\nicefrac{1}2}_h|U^1_h|^2 \dif x \right) \bigg |,
\end{aligned}
\end{equation}}
Table \ref{invsVC} shows the corresponding conservation of mass \eqref{DMC} and balance of energy \eqref{localCE}.  We observed that both are recovered to double precision of accuracy for all values of $\ep$. 
\begin{table}[htp]
\caption{Errors in the conservation laws : variable $p(t), \ q(t)$}\vspace{-5mm}
\begin{center}
\begin{tabular}{|c|c|c||c|c||c|c||} \hline
   &  \multicolumn{2}{c||}{$\varepsilon=1$} & \multicolumn{2}{c||}{$\varepsilon=0.1$}  & \multicolumn{2}{c||}{$\varepsilon=0.01$}  \\ \hline
$t_n$  & 	$ \mathcal{M}_e^n$    & 	$ \E_{e,loc}^n$   & 	$ \mathcal{M}_e^n$     & 	$ \E_{e,loc}^n$  &  $ \mathcal{M}_e^n$    &  $ \E_{e,loc}^n$ \\ \hline 
0	&	2.05e-15	&	6.69e-16	&		4.49e-15	&	6.27e-15	&		6.55e-15 	&	1.08e-15	\\ \hline 
1	&	1.91e-14	&	5.81e-16	& 		9.81e-15	&	3.29e-15	&		5.55e-15	&	1.74e-14	\\ \hline
2	&	5.73e-14	&	3.62e-16	&		1.97e-14	&	4.02e-15	&		8.91e-15	&	2.86e-14	\\ \hline
3	&	9.38e-14	&	2.71e-16	&		1.47e-13	&	4.39e-15	&		1.78e-14	&	4.39e-14	\\ \hline
\end{tabular}
\end{center}
\label{invsVC}%\vspace{-10mm}
\end{table}
%end RED
%

{\color{black}
\subsection{Variable time-step $k_n$}
We examine now the effect of variable time-step in the mass conservation \eqref{cl1} and energy balance \eqref{cl2} of system \eqref{SP}.  We take \eqref{CLIC} as an initial condition, $\alpha=\beta=5$ and $\epsilon=0.01$. The domain $\varOmega = (-1,1)^2$ is discretized by a uniform grid $\Th$ consisting of squares with sides of length $h=0.015625$, and we consider cubic finite elements on $\Th$ resulting a spatial discretization error which is almost negligible. We take  $\tau=0$ and  final time  $T = 3$.  We split the time interval $[0,T] = \cup_{j=1}^{12} [T_{j-1}, T_{j}]$ with $T_j = j/4, \ 0\le j \le 12$. In each subinterval $[T_{j-1}, T_{j})$ we use a different time step $k_j = j\cdot 1.25e-03, \  1\le j \le 12$. 
We monitor the error in discrete mass conservation \eqref{DMC} by means of $\mathcal{M}_e^j$ which was defined in \eqref{dMC}. Moreover we monitor error in the energy balance law \eqref{DEC} by means of the global error $\E_{e,gl}^j$ defined in \eqref{dEC}, and of the  local error $\mathcal{E}^{j}_{e,loc}$ defined in \eqref{dECloc}.

At this point one should also recall that, according to Proposition \ref{lm:jump}, the size of $ \E_{e,loc}^{j}$ at the points of change of time-step size is expected to be equal to the residual
\begin{equation}
\mathcal{R}^{j}:=\bigg|\frac{q(t_{j-1/2})}{4}\frac{k_{j-1}-k_{j}}{k_{j-1}+k_{j}}k_{j}^2\|\bar\p\nabla V_h^{j-1/2}\|^2\bigg|.
\end{equation}

%%
%\begin{table}[htp]
%\caption{Errors in the conservation laws : variable time-step }
%\begin{center}
%%\small
%\begin{tabular}{|c|c|c|c|c|c|} \hline
%$T_n$                 &    $k_n$         & $ \mathcal{M}_e^n$ & $ \E_{e,gl}^{n+1}$  &  $  \E_{e,loc}^{n+1}$  &  $\mathcal{R}^{n+1}$  \\ \hline
%$0.00 - 0.25$    & 	1.250e-03	& 	2.99e-15	& 	3.11e-15	& 	2.88e-11		& 	 1.28e-11	\\ \hline
% $0.25 - 0.50$	& 	2.500e-03	& 	7.54e-14	& 	2.93e-11	& 	1.92e-10	    & 	1.22e-10	\\ \hline
% $0.50 - 0.75$	& 	3.750e-03	& 	9.30e-14	& 	2.21e-10	& 	6.09e-10& 	4.42e-10	\\ \hline
%$0.75 - 1.00$ & 	5.000e-03	& 	1.07e-13	& 	8.31e-10	&   1.40e-09	& 	1.09e-09	\\ \hline
%$1.00 - 1.25$ & 	6.250e-03	& 	1.23e-13	& 	2.23e-09	& 	2.69e-09        & 	2.20e-09	\\ \hline
%$1.25 - 1.50$ 	& 	7.500e-03	& 	1.28e-13	& 	4.92e-09	& 	4.66e-09		& 	3.92e-09	\\ \hline
%$1.50 - 1.75$ 	& 	8.750e-03	& 	1.39e-13	& 	9.59e-09	& 	7.37e-09	& 	6.32e-09	\\ \hline
%$1.75 - 2.00$ 	& 	1.000e-02	& 	1.41e-13	& 	1.70e-08	& 	1.10e-08	& 	9.71e-09	\\ \hline
%$2.00 - 2.25$ 	& 	1.125e-02	& 	1.49e-13	& 	2.80e-08	& 	1.61e-08	& 	1.42e-08	\\ \hline
%$2.25 - 2.50$ 	& 	1.250e-02	& 	1.53e-13	& 	4.40e-08	& 	2.24e-08	& 	2.00e-08	\\ \hline
%$2.50 - 2.75$ 	& 	1.375e-02	& 	1.59e-13	& 	6.64e-08	& 	3.04e-08	& 	2.72e-08	\\ \hline
%$2.75 - 3.00$ 	& 	1.500e-02	& 	1.63e-13	& 	9.68e-08	& 	-- & -- 		\\ \hline
%\end{tabular}
%\end{center}
%\label{invsVT}
%\end{table}

%
\begin{table}[htp]
\caption{Errors in the conservation laws : variable time-step }
\begin{center}
%\small
\begin{tabular}{|c|c|c|c|c|c|c|} \hline
$T_{j}$                 &    $k_{j}$ & $k_{j+1}$         & $ \mathcal{M}_e^{j}$ & $ \E_{e,gl}^{j}$  &  $  \E_{e,loc}^{j}$  &  $\mathcal{R}^{j}$  \\ \hline
$0.25$    & 	1.250e-03	&2.500e-03	& 	2.99e-15	& 	3.11e-15	& 	2.88e-11		& 	 1.28e-11	\\ \hline
 $0.50$	& 	2.500e-03	&3.750e-03	& 	7.54e-14	& 	2.93e-11	& 	1.92e-10	    & 	1.22e-10	\\ \hline
 $0.75$	& 	3.750e-03	&5.000e-03	& 	9.30e-14	& 	2.21e-10	& 	6.09e-10& 	4.42e-10	\\ \hline
$1.00$ & 	5.000e-03	&6.250e-03	& 	1.07e-13	& 	8.31e-10	&   1.40e-09	& 	1.09e-09	\\ \hline
$1.25$ & 	6.250e-03	&7.500e-03	& 	1.23e-13	& 	2.23e-09	& 	2.69e-09        & 	2.20e-09	\\ \hline
$1.50$ 	& 	7.500e-03	&8.750e-03	& 	1.28e-13	& 	4.92e-09	& 	4.66e-09		& 	3.92e-09	\\ \hline
$1.75$ 	& 	8.750e-03	&1.000e-02	& 	1.39e-13	& 	9.59e-09	& 	7.37e-09	& 	6.32e-09	\\ \hline
$2.00$ 	& 	1.000e-02	&1.125e-02	& 	1.41e-13	& 	1.70e-08	& 	1.10e-08	& 	9.71e-09	\\ \hline
$2.25$ 	& 	1.125e-02	&1.250e-02	& 	1.49e-13	& 	2.80e-08	& 	1.61e-08	& 	1.42e-08	\\ \hline
$2.50$ 	& 	1.250e-02	&1.375e-02	& 	1.53e-13	& 	4.40e-08	& 	2.24e-08	& 	2.00e-08	\\ \hline
$2.75$ 	& 	1.375e-02	&1.500e-02	& 	1.59e-13	& 	6.64e-08	& 	3.04e-08	& 	2.72e-08	\\ \hline
$3.00$ 	& 	1.500e-02	& -- & 	1.63e-13	& 	9.68e-08	& 	-- & -- 		\\ \hline
\end{tabular}
\end{center}
\label{invsVT}
\end{table}
Our findings are presented in Table \ref{invsVT}. Each row $j , \ 1\le j\le 12$ corresponds to the time interval $[T_{j-1}, T_j)$ and $T_{j}$ is printed in the first column, while the corresponding time-step $k_{j}$ is shown in the second column. The time-step size changes right after $T_{j};$ for the computation of the numerical solution in $[T_j, T_{j+1})$ with new time-step $k_{j+1}$ appears (cf. equation \eqref{CNrelax}); this $k_{j+1}$ is printed in the third column. The global mass error at  $T_{j}$ is printed in the fourth column, and the global energy error at $T_{j}$ is printed in the fifth column. The local error in energy due to the change of time-step size around $T_{j}$ is printed in the sixth column. The predicted residual, which is expected to be equal to the local energy error, is printed in the last column.  We observe that conservation of mass is essentially unaffected by the time-step change while the energy loss in monitored accurately by the residual derived in Proposition \ref{lm:jump}.  We note that, during this computation the time-step changed 11 times and grew by a factor of 20. 
}
%\BLUE{
%\begin{remark}
%An interested numerical observation from table 5 is the following: for long time $\mathcal{R}^{j}$ seems to behave like 
%$$
%\mathcal{R}^{j} \sim \frac{1}{11} k_j^{3.5} 
%$$ 
%where $11$ is the number of times we changed the time-step !!! I don't think it's a coincidence.
%\end{remark}
%}

\subsection{A cosmological example} \label{sec:cosmoex}
\begin{figure}
		\includegraphics[width=8.5cm, height=8cm]{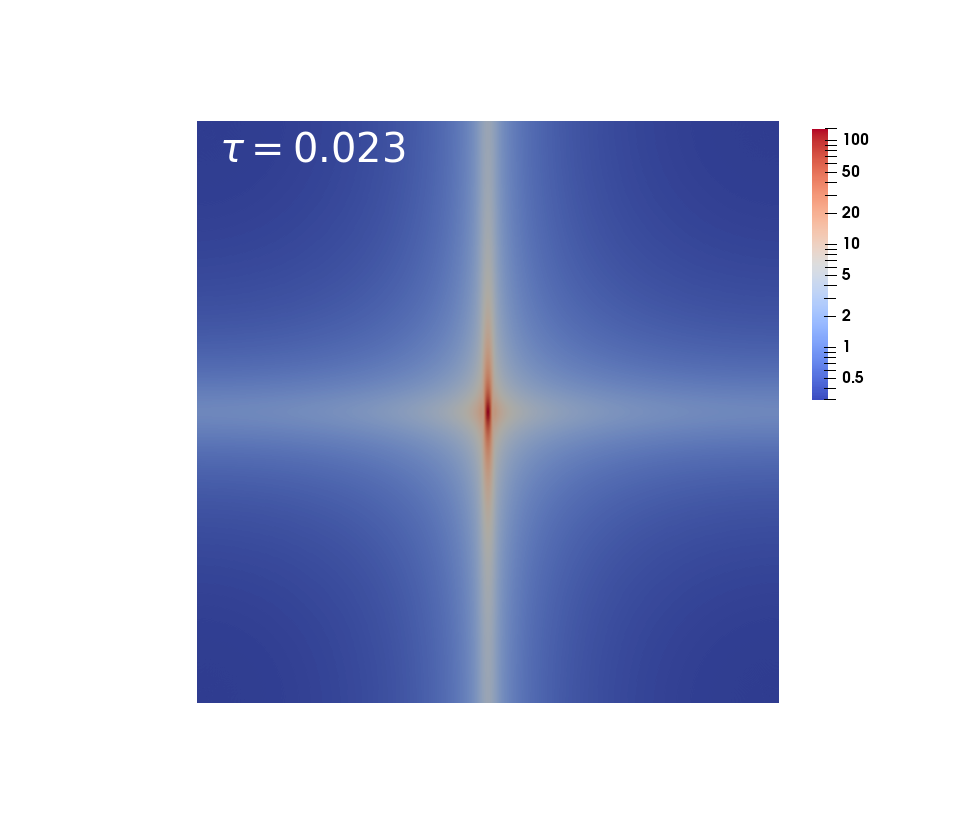} \hspace{-2.25cm}
		\includegraphics[width=8.5cm, height=8cm]{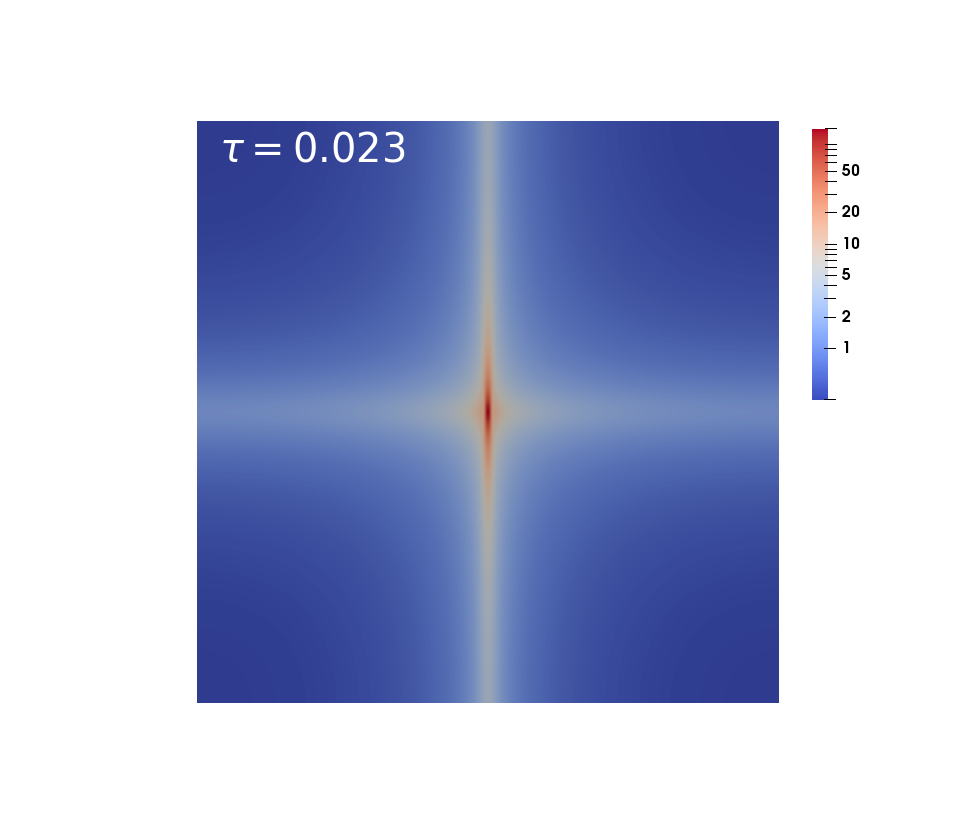} 
		
		\vspace{-1.5cm}
		
		\includegraphics[width=8.5cm, height=8cm]{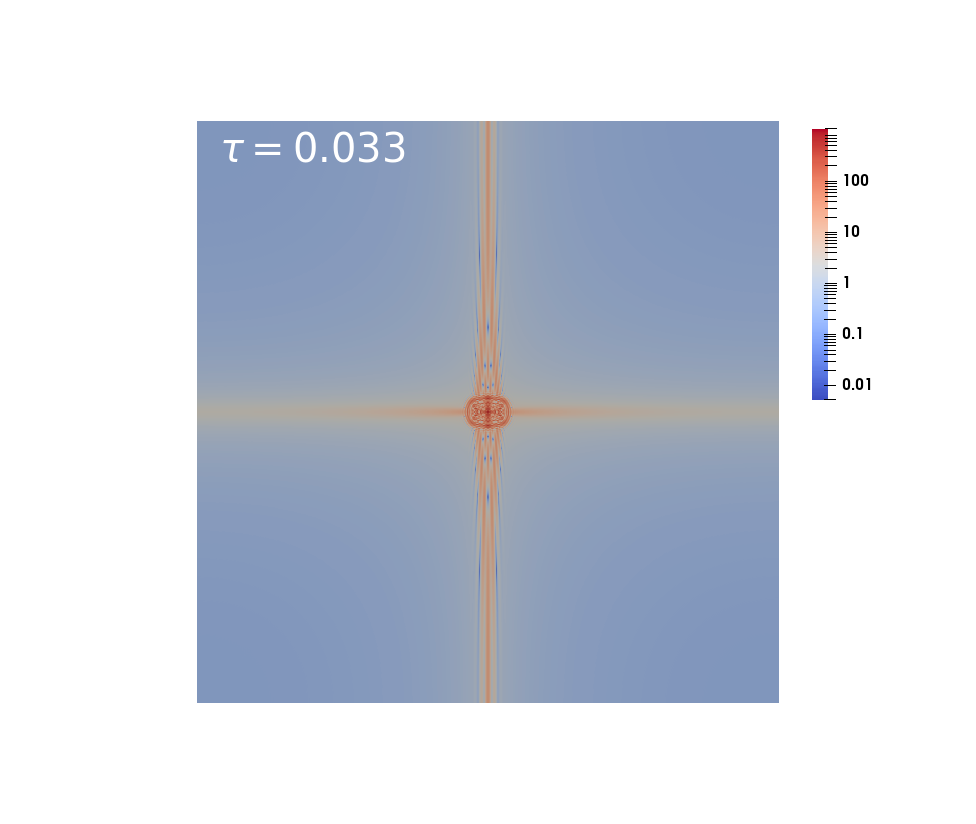} \hspace{-2.25cm}
		\includegraphics[width=8.5cm, height=8cm]{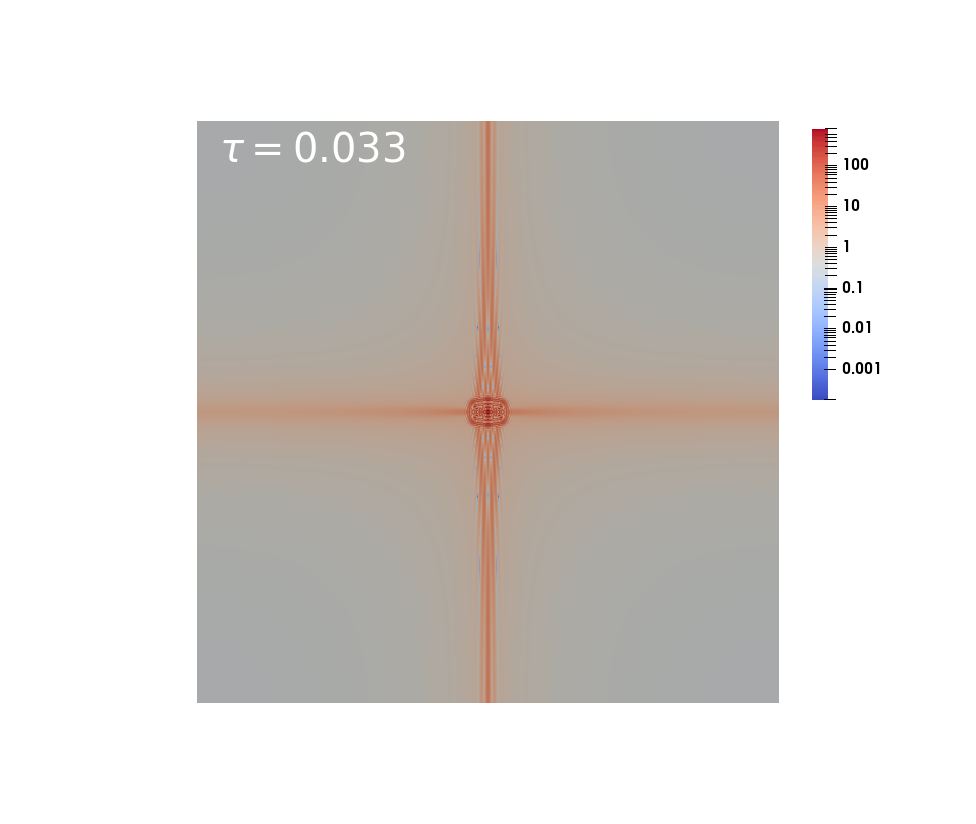}
		
		\vspace{-1.5cm}
		
		\includegraphics[width=8.5cm, height=8cm]{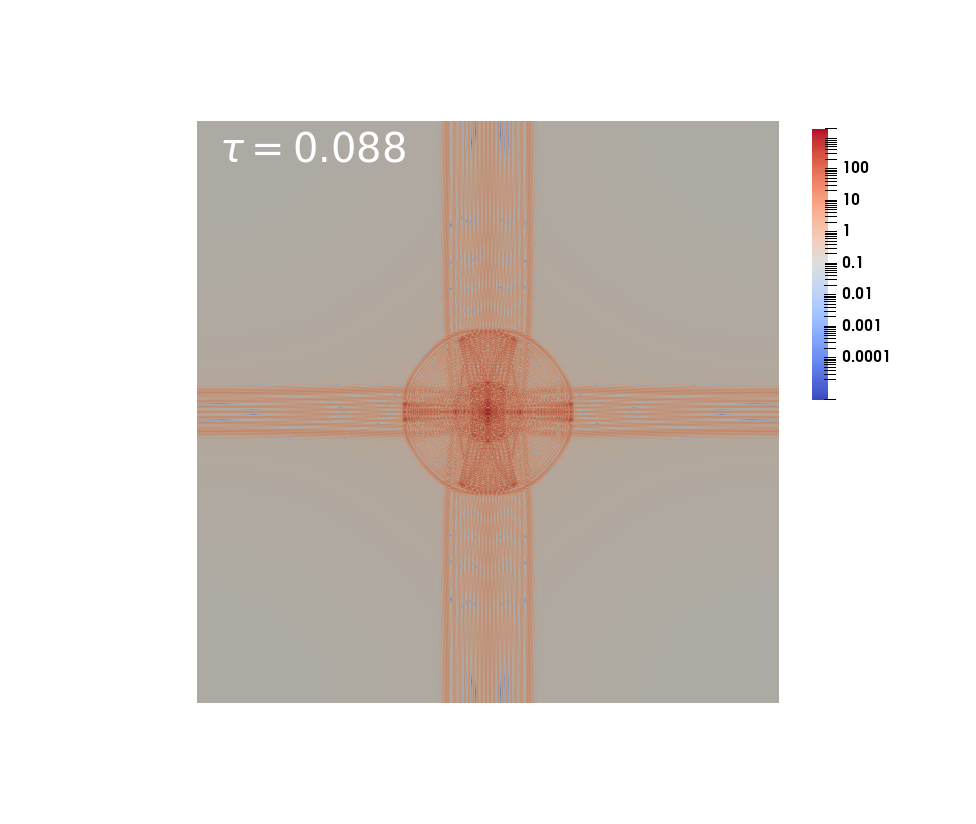} \hspace{-2.25cm}
		\includegraphics[width=8.5cm, height=8cm]{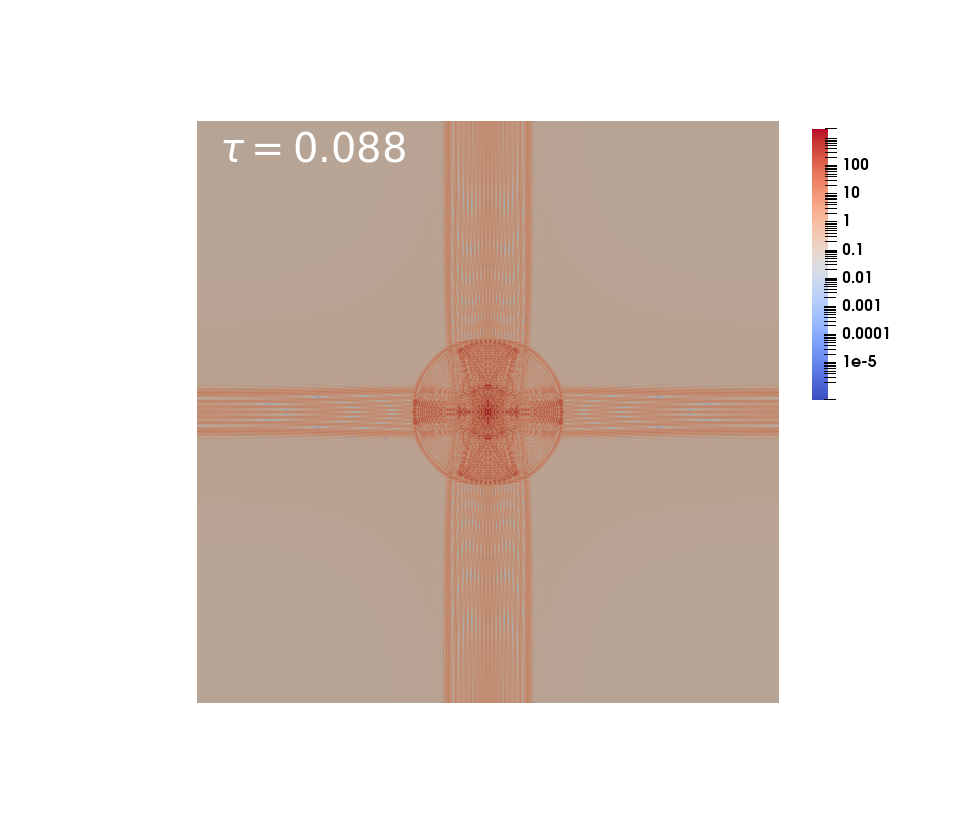} \vspace{-1cm}
	\caption{The density $|U_h(t)|^2$ (logarithmic scale) at $t=0.023, \ 0.033, \ 0.088$ : $1024\!\times\! 1024$ grid (left), $2048\!\times\!2048$ grid (right).}
\label{UhMod}	
\end{figure}

One application of the Schr\"odinger-Poisson system \eqref{GeSP} comes from the field of cosmology. Indeed, the $d$-dimensional Schr\"odinger-Poisson system \eqref{GeSP} can be used as an approximation to the computationally expensive $2d$-dimensional Vlasov-Poisson system used to describe collisionless self-gravitating matter \cite{KVS}.
\begin{figure}[h]%
	\includegraphics[width=0.49\textwidth]{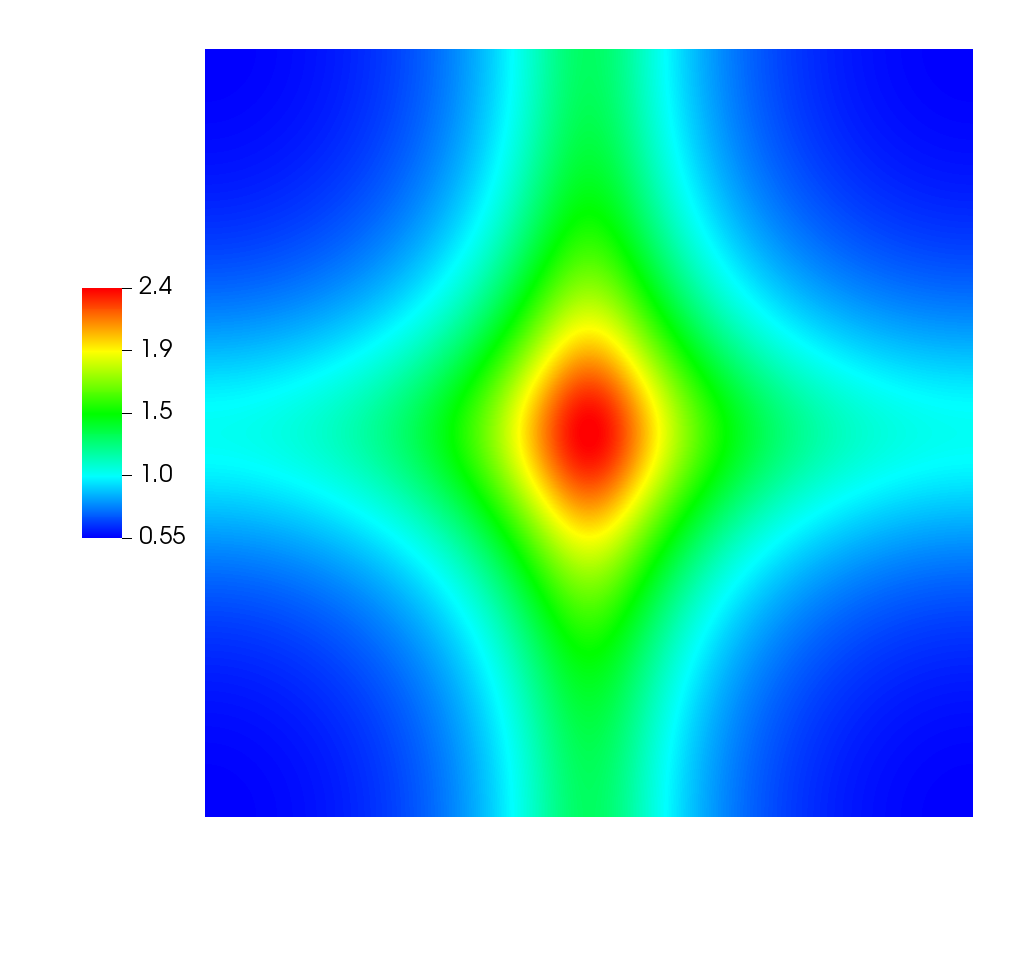}
	\includegraphics[width=0.49\textwidth]{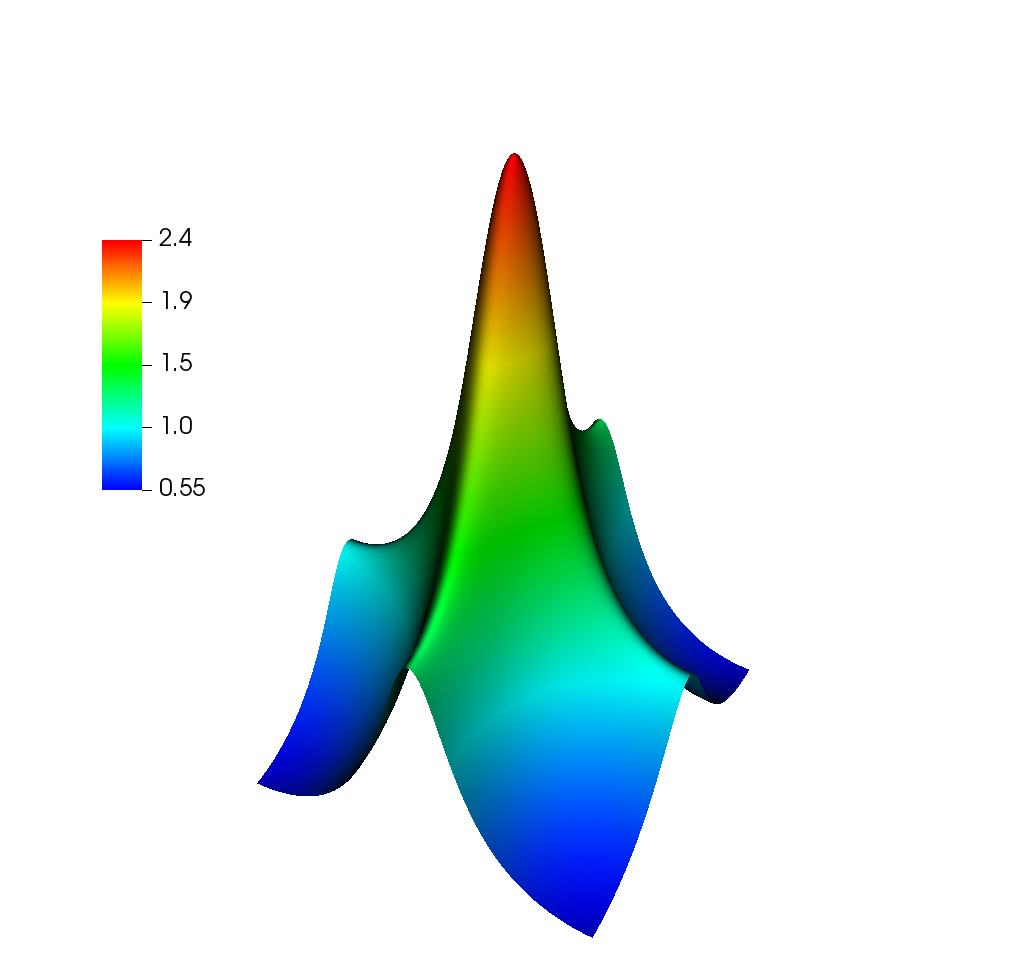}
	\caption{The initial density $|u_0|^2$ used in \eqref{SPnew}.} \label{u0Mod}
\end{figure}%
\noindent For this numerical experiment, we proceed as in \cite{KVS}. Firstly, we take $\mu=\|u_0\|^2=1$ and periodic boundary conditions in \eqref{GeSP}.  The initial time is $\tau=0.01$ and the final time $T=0.088.$ Time dependent coefficients are used, namely $p(t) = \frac{\ep}{2 t^{\nicefrac{3}{2}}}$, $q(t) = \frac{\beta}{\ep t^{\nicefrac{1}{2}}},$ where $\beta=1.5$ and $\varepsilon = 6e$\,--\,$5.$ Thus our cosmological example reads as
\begin{equation}
	\label{SPnew} \left \{
	\begin{aligned}
		& u_{t}-\frac{\ii\ep}{2t^{\nicefrac32}}\Delta u +\frac{\ii\beta}{\ep t^{\nicefrac{1}{2}}} \V u=0,    &&\qquad\mbox{in ${\varOmega}\!\times\! (\tau,T)$,}& \\
		&\Delta\V = |u|^2-1,  &&\qquad\mbox{in ${\varOmega}\!\times\! (\tau,T)$,}&\\
		& u(\xv,0) = u_0(\xv),  && \qquad\mbox{in $\varOmega$,} \\
		&   \|u_0\|_{L^2}^2=1 :  u, \V \text{ periodic } && \qquad \mbox{on $\partial \varOmega\!\times\! (\tau,T)$.}&  \\
	\end{aligned}
	\right.
\end{equation}%
This problem is a special case of system \eqref{GeSP}, and therefore satisfies the mass conservation  \eqref{cl1} and the energy balance \eqref{cl2}, which for the particular $p(t)$ and $q(t)$ takes the form:
\begin{equation*}
 \frac{\ep}{2 t^{\nicefrac{3}{2}}} \frac{\dif}{\dif t} \E_k(t) + \frac{\beta}{\ep t^{\nicefrac{1}{2}}} \frac{\dif}{\dif t} \left(2 \E_{\V_2}(t) -  \E_{\V_1}(t)\right) = 0 \implies \frac{\ep^2}{2} \frac{\dif}{\dif t} \E_k(t) + t \beta \frac{\dif}{\dif t}\left( 2\E_{\V_2}(t) -  \E_{\V_1}(t)\right) = 0 .
\end{equation*}
The fully discrete relaxation scheme \eqref{CNrelax} with constant time-step $k_n = k, \ \forall n$ for \eqref{SPnew} is 
\begin{equation}
\label{CNrelaxC} \left \{
\begin{aligned}
&\overline{\partial} U_h^n-\frac{\ii\varepsilon}{2t_{n-\nicefrac{1}{2}}^{\nicefrac32}} \Delta_h^n U_h^{n-\nicefrac{1}{2}}+\frac{\ii\beta}{\ep t_{n-\nicefrac12}^{\nicefrac12}}\Pc_h^{n}\left(V_h^{n-\nicefrac{1}{2}}U_h^{n-\nicefrac{1}{2}}\right)=0,\\
& \Delta_h^n V_h^{n-\nicefrac{1}{2}} = \Phi_h^{n-\nicefrac{1}{2}}, \\
&\Phi_h^{n-\nicefrac{1}{2}} = 2 \Pc_h^{n} (|U_h^{n-1}|^{2}-1) -  \Phi_h^{n-\nicefrac{3}{2}}.
\end{aligned}
\right.
\end{equation}
Scheme \eqref{CNrelaxC} satisfies the discrete mass conservation \eqref{DMC} and the local discrete energy balance  \eqref{localCE}, which for the particular $p(t)$ and $q(t)$ takes the form:
\begin{equation}
\ep^2 \overline{\partial}(\|\nabla U_h^n\|^2) + \beta t_{n-\nicefrac{1}{2}} \left(   \overline{\partial}(\|\nabla V_h^{n-\nicefrac{1}{2}}\|^2) + 2\int_{\varOmega}  \overline{\partial}(V_h^{n - \nicefrac{1}{2}}(|U_h^{n}|^2-\mu)) \dif x  \right) = 0, \forall n .
\end{equation}
We now consider the benchmark case of a \emph{sine wave collapse}.  The initial condition is as given in \cite{KVS} and its position density is displayed in Figure \ref{u0Mod}.
The domain $\varOmega$ is discretized with linear finite elements over one of two different uniform grids: a $1024\!\times\! 1024$ and a $2048\!\times\! 2048$ grid while the time domain $(\tau, T)$ is discretized using 1560 uniform time steps yielding a time step size of $k =  5e$\,--\,$5$. Results of the numerical simulations are shown in Figure \ref{UhMod}.  For comparison purposes, we plot the density $|U_h(t_n)|^2$ at three different time instances ($t_n=0.023, \ 0.033, \ 0.088$) all of which are in {excellent} agreement with the plots in \cite{KVS}. 

\begin{figure}[h]
		\includegraphics[scale=0.44]{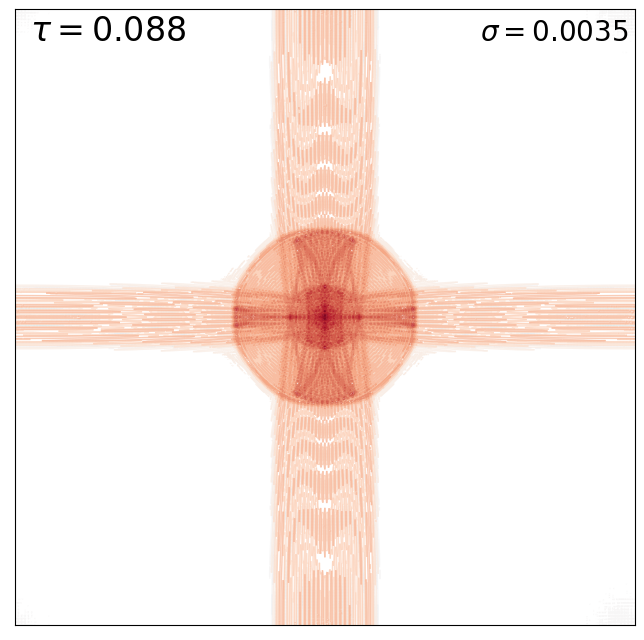}  \hspace{0.5cm}
		\includegraphics[scale=0.44]{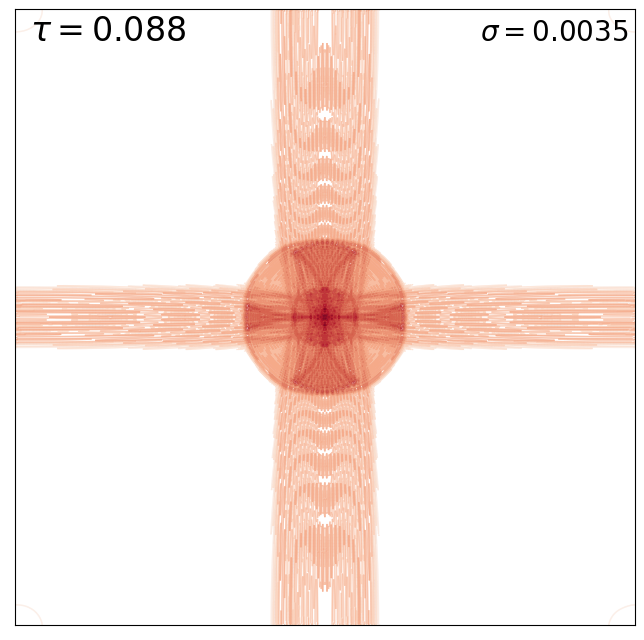}
	\caption{The density $|U_h(t)|^2$ (logarithmic scale) at $t=T$  with Gaussian filtering of width $\sigma=0.0035$ : $1024\!\times\! 1024$ grid (left), $2048\!\times\! 2048$ grid (right).} \label{GFUhMod}
\end{figure}

Due to the small value of $\varepsilon$, the wavefunction is highly oscilliatory which can be readily seen in the $1024\!\times\!1024$ grid but is much more apparent in the $2048\!\times\!2048$ grid. We thus postprocess the density by applying a Gaussian filter of width $\sigma = 0.0035$ which eliminates spurious oscillations  {\bf --} this is shown in  Figure \ref{GFUhMod} at the final time $\tau=T$ for both grids.

%
%\begin{figure}
%		\includegraphics[width=8.5cm, height=8cm]{Fig_e6_5_a0_023_1024.png} \hspace{-2.25cm}
%		\includegraphics[width=8.5cm, height=8cm]{Fig_e6_5_a0_023_2048.png} 
%		
%		\vspace{-1.5cm}
%		
%		\includegraphics[width=8.5cm, height=8cm]{Fig_e6_5_a0_033_1024.png} \hspace{-2.25cm}
%		\includegraphics[width=8.5cm, height=8cm]{Fig_e6_5_a0_033_2048.png}
%		
%		\vspace{-1.5cm}
%		
%		\includegraphics[width=8.5cm, height=8cm]{Fig_e6_5_a0_088_1024.png} \hspace{-2.25cm}
%		\includegraphics[width=8.5cm, height=8cm]{Fig_e6_5_a0_088_2048.png} \vspace{-1cm}
%	\caption{The density $|U_h(t)|^2$ (logarithmic scale) at $t=0.023, \ 0.033, \ 0.088$ : $1024\!\times\! 1024$ grid (left), $2048\!\times\!2048$ grid (right).}
%\label{UhMod}	
%\end{figure}%
%%

For this example, the density is conserved up to double precision while the energy is conserved to at least  $8$ digits of accuracy, see Table \ref{MCeEBe}. The loss of accuracy in the energy balance is mainly due to the loss of periodicity for  $\nabla \V$ at the discrete level. Indeed, the periodicity for wavefuction and potential is preserved at the discrete level, however it is lost for  $\nabla V_h^n$ which is a  crucial assumption in \eqref{ce2} for proving \eqref{localBE}.  We have computed the $L^2$-norm of the difference between the values of the $\nabla V_h^n, \ n=1,\dots,$ along the corresponding horizontal and vertical boundaries of the domain and it is found to vary from $10^{-7} - 10^{-10}$ depending on the grid size.

\begin{table}[htbp]
 \caption{Mass conservation error(MCe) and Energy balance error (EBe)}
  \begin{center}
    \begin{tabular}{|c||c|c||c|c||c|c||}\hline
    Grid &  \multicolumn{2}{c||}{$512\times512$} & \multicolumn{2}{c||}{$1024\times1024$}  & \multicolumn{2}{c||}{$2048\times2048$}  \\ \hline
    t              &  MCe         & EBe             &  MCe          & EBe            &  MCe           & EBe \\ \hline\hline
    $0.023$  & 4.329e-15  &  1.222e-08  & 1.010e-13  &  9.039e-10  &  1.852e-13  &  1.131e-11\\
    $0.033$  & 1.643e-14  &  2.427e-08  & 1.210e-13  &  9.675e-10  &  1.912e-13  &  2.323e-11\\
    $0.088$  & 3.775e-14  &  5.065e-08  & 1.386e-13  &  9.273e-10  &  5.332e-13  &  6.037e-11 \\
    \hline
    \end{tabular}\\[2ex]
  \end{center}
  \label{MCeEBe}%
\end{table}%}
%end RED
%

%{\color{black}
%\subsection{Energy balance remainder for variable time-steps} \label{sec:numjump}
%
%In this section we examine by how much the energy balance is off at the discrete level when the time-step size changes. Indeed, according to Lemma \ref{sec:conslnum}
%
%Here we use the same problem as in Section \ref{sec:order}
%
%}

%\section{Conclusions}
%
%We proposed a new numerical method \eqref{RCNG} for the numerical solution of  a class of Schr\"odingerÃ¢ÂÂPoisson type systems \eqref{GeSP} which is order two in time and of optimal order in space. The numerical method is linearly implicit, so no nonlinear iteration (e.g. Newton's method) is required, making the method fast as well as accurate. The numerical method conserves mass and satisfies the same energy balance law as the continuous problem.   
%
\section*{Acknowledgements}
The authors acknowledge the support from the Carnegie Trust Research Incentive Grant RIG008215. I.K. would also like to acknowledge the support from London Mathematical Society through an Emmy Noether Fellowship. 
In addition, Th.\,K. and I.K. thank the Edinburgh Mathematical Society for the Covid Recovery Fund that allowed for the completion and the submission of this paper.
 Moreover, the authors would like to express their gratitude to Dr. K. Vattis, Prof. C. Skordis and especially Dr. M. Kopp for their valuable help and support in setting up the cosmological example reported in Section 5. Finally, the authors would like to thank the anonymous reviewers for their valuable comments and suggestions.

%%%%%%%%%%%%%%%%%%%%%%%%%%%
%%% BIBLIOGRAPHY
%%%%%%%%%%%%%%%%%%%%%%%%%%%

%
%-------------------------------------------------------------------------------------------------------------------------

\begin{thebibliography}{10}
	 
%
\bibitem{AMP}
N.B. Abdallah, F. M\'ehats, O. Pinaud, \emph{On an open transient Schr\"odinger-Poisson system}, 
Math. Models Methods Appl. Sci.  \textbf{15},  667--688, 2005.
%
\bibitem{ADK}
{\color{black}
G.\ Akrivis, V.\ Dougalis, O.\ Karakashian, \emph{On fully
discrete Galerkin methods of second-order temporal accuracy for
the nonlinear Schr\"{o}dinger equation},
\newblock Numer.\ Math.\ \textbf{59}, 31-53, 1991.
}
%
\bibitem{ADKM}
{\color{black}
G.D.\ Akrivis, V.A.\ Dougalis, O.A.\ Karakashian, W.R.\
 McKinney, \emph{Numerical approximation of blow-up of radially
symmetric solutions of the nonlinear Schr\"{o}dinger equation},
\newblock SIAM J.\ Sci.\ Comput.\ \textbf{25}, 186--212, 2003.
}
%
\bibitem{AkrivisLi}
{\color{black} G.\ Akrivis, D.\ Li, \emph{Structure-preserving Gauss methods for the nonlinear Schr\"odinger equation}, Calcolo \textbf{58}, 1--25, 2021.
}
\bibitem{AS}
E. Arriola, J. Soler, \emph{A variational approach to the Schr\"odinger-Poisson System: Asymptotic behaviour, breathers, and stability}, J. Stat. Phys. \textbf{103}, 1069--1106, 2001.
%
\bibitem{AKKT}
W. Auzinger, T. Kassebacher, O. Koch, M. Thalhammer, \emph{Convergence of a Strang splitting finite element discretization for the Schr\"odinger-Poisson equation},  ESAIM Math. Model. Numer. Anal. \textbf{51}, 1245--1278, 2017.
%
\bibitem{Bao}
{\color{black}
W. Bao, D. Jaksch, P.A. Markowich,   \emph{Numerical solution of the Gross-Pitaevskii equation for Bose-Einstein condensation}, J. Comput. Phys. \textbf{187}, 318--342, 2003.
}
\bibitem{BHK}
W. Bangerth, R. Hartmann, G. Kanschat, \emph{{\tt deal.II} -- A general-purpose object-oriented finite element library}, ACM Transactions on Mathematical Software, 33(4), article 24, 2007.
%
\bibitem{BNS}
W. Bao, N. Mauser,  H.P. Stimming, \emph{Effective one particle quantum dynamics of electrons: A numerical study of the Schr\"odinger-Poisson-X$\alpha$ model}, Comm. Math. Sciences \textbf{1}, 809--828, 2003.
%
\bibitem{Besse}
Ch. Besse, \emph{A relaxation scheme for the nonlinear Schr\"odinger equation},  SIAM J. Numer. Anal. \textbf{42}, 934--952, 2004.
%
\bibitem{Besse2}
Ch. Besse, S. Descombes, G. Dujardin, I. Lacroix-Violet, \emph{Energy-preserving methods for nonlinear Schr\"odinger equations}, IMA J. Numer. Anal. \textbf{41}, 618--653, 2021.
%
\bibitem{Besse3}
{\color{black}
C.\ Besse, G.\ Dujardin, I. Lacroix-Violet,   \emph{High order exponential integrators for nonlinear Schr\"odinger equations with application to rotating Bose--Einstein condensates}, SIAM J. Numer. Anal. \textbf{55}, 1387--1411, 2017.
}
%
 
%
%
\bibitem{BEGMY}
C. Bardos, L. Erd\"os, F. Golse, N. Mauser, H-T Yau, \emph{Derivation of the Schr\"odinger-Poisson equation from the quantum N-body problem}, C. R. Acad. Sci. Paris, Ser. I \textbf{334}, 515-520, 2002.
%
\bibitem{Berland}
{\color{black}
H. Berland, A.L. Islas, C.M. Schober,   \emph{Conservation of phase space properties using exponential integrators on the cubic Schr\"odinger equation}, J. Comput. Phys. \textbf{225}, 284--299, 2007.
}
\bibitem{BTGIFG}
P. Bertrand, N. Van Tuan, M. Gros, B. Izrar, M. Feix,  J. Gutierrez, \emph{Classical Vlasov plasma description through quantum numerical methods}, J. Plasma Phys. \textbf{23}, 401-422, 1980.

%
\bibitem{BILZ}
S. Bohun, R. Illner, H. Lange, P.F Zweifel, \emph{Error estimates for Galerkin approximations to the periodic Schr\"odinger-Poisson system}, ZAMM Journal of applied mathematics and mechanics/ Zeitschrift f\"ur angewandte Mathematik und Mechanik \textbf{76},  7--13, 1996
%
\bibitem{BM}
F. Brezzi,  P.A. Markowich, \emph{The three-dimensional Wigner-Poisson problem: Existence, uniqueness and approximation}, Math. Methods Appl. Sci. \textbf{14}, 35--61, 1991
%

%
 
%
\bibitem{Castella}
F. Castella, \emph{$L^2$-solutions to the Schr\"odinger-Poisson system: existence, uniqueness, time behaviour and smoothing effects}, Math. Mod. Meth. Appl. Sci. \textbf{7}, 1051-1083, 1997.
%
\bibitem{Cazenave}
T. Cazenave, \emph{Semilinear Schrodinger Equations} (Vol. 10). American Mathematical Soc. 2003.

\bibitem{Chartier}
{\color{black}
P. Chartier, N.J. Mauser, F. M\'ehats,   Y. Zhang, Y., \emph{Solving highly-oscillatory NLS with SAM: numerical efficiency and long-time behavior}, Discrete Contin. Dyn. Syst.-S \textbf{9}, 1327, 2016.
}

\bibitem{DW}
G. Davies, L. Widrow, \emph{Test-bed simulations of collisionless, self-gravitating systems using the Schr\"odinger method}, The Astrophysical Journal \textbf{485}, 484, 1997.
\bibitem{Delfour}
{\color{black}
M. Delfour, M. Fortin, G. Payr, \emph{Finite-difference solutions of a non-linear Schr\"odinger equation}, J. Comput. Phys. \textbf{44}, 277--288, 1981.
}
\bibitem{Deg1}
{\color{black}
M. Dehghan, V. Mohammadi,  \emph{A numerical scheme based on radial basis function finite difference (RBF-FD) technique for solving the high-dimensional nonlinear Schr\"odinger equations using an explicit time discretization: Runge-Kutta method}, Comput. Phys. Commun. \textbf{217}, 23--34, 2017.
}
\bibitem{Deg2}
{\color{black}
M. Dehghan, A. Taleei, \emph{A compact split-step finite difference method for solving the nonlinear Schr\"odinger equations with constant and variable coefficients}, Comput. Phys. Commun. \textbf{181}, 43--51, 2010.
}
%
\bibitem{EZ}
M. Ehrhardt, A. Zisowsky, \emph{Fast calculation of energy and mass preserving solutions of Schr\"odinger-Poisson systems on unbounded domains}, J. Comput. Appl. Math. \textbf{187}, 1-28, 2006.

\bibitem{Fei}
{\color{black}
Z. Fei, V.M. P\'erez-Garc\'i­a, L. V\'azquez,  \emph{Numerical simulation of nonlinear Schr\"odinger systems: a new conservative scheme} Appl. Math. Comput. \textbf{71}, 165--177, 1995.
}
%
\bibitem{Hederi}
{\color{black}
M. Hederi, A.L. Islas, K. Reger,  C.M. Schober,  \emph{Efficiency of exponential time differencing schemes for nonlinear Schr\"odinger equations}, Math. Comput. Simul. \textbf{127}, 101--113, 2016.
}
\bibitem{Deg3}
{\color{black}
M. Ilati, M. Dehghan,  \emph{DMLPG method for numerical simulation of soliton collisions in multi-dimensional coupled damped nonlinear Schr\"odinger system which arises from Bose-Einstein condensates}, Appl. Math. Comput. \textbf{346}, 244--253, 2019.
}
\bibitem{IZL}
R. Illner, P.F. Zweifel, H. Lange, \emph{Global existence, uniqueness and asymptotic behaviour of solutions of the Wigner-Poisson and Schr\"odinger-Poisson systems}, Math. Meth. Appl. Sci. \textbf{17}, 349--376, 1994.
%
\bibitem{JWY}
S. Jin , H. Wu, X. Yang, \emph{A numerical study of the Gaussian beam methods for Schr\"odinger-Poisson
equations}, J. Comput. Appl. Math. \textbf{28}, 261-272, 2010.
%

\bibitem{Karner}
M. Karner, A. Gehring, S. Holzer, M. Pourfath, M. Wagner, W. Goes, M. Vasicek, O. Baumgartner, C. Kernstock, K. Schnass, G. Zeiler,   \emph{A multi-purpose Schr\"odinger-Poisson solver for TCAD applications} J. Comput. Electron. \textbf{6}, 179--182, 2007.

\bibitem{Kara}
{\color{black}
O. Karakashian,  Ch. Makridakis, \emph{A space-time finite element method for the nonlinear Schr\"odinger equation: the discontinuous Galerkin method}, Math. Comp. \textbf{67}, 479--499, 1998.
}

\bibitem{KM2}
{\color{black} 
O. Karakashian, Ch. Makridakis, \emph{A space-time finite element method for the nonlinear
Schr\"odinger equation: the continuous Galerkin method}, SIAM J. Numer. Anal. \textbf{36},
1779--1807, 1999.}

\bibitem{KK}
Th. Katsaounis, I. Kyza, \emph{A posteriori error analysis for evolution nonlinear Schr\"odinger equations up to the critical exponent}, SIAM J. Numer. Anal. \textbf{56}, 1405--1434, 2018.
%
\bibitem{KVS}
M. Kopp, K. Vattis, C. Skordis, \emph{Solving the Vlasov equation in two spatial dimensions with the Schr\"odinger method}, Phys. Rev. D \textbf{96}, 123532, 2017.
%
\bibitem{LC}
T. Lu, W. Cai, \emph{A Fourier spectral-discontinuous Galerkin method for time-dependent 3-D Schr\"odinger-Poisson equations with discontinuous potentials}, J. Comput. Appl. Math. \textbf{220}, 588-614, 2008.
%
\bibitem{Lubich}
C. Lubich, \emph{On splitting methods for Schr\"odinger-Poisson and cubic nonlinear Schr\"odinger equations}, Math. Comp. \textbf{77}, 2141--2153, 2008.
%

%
\bibitem{MRS}
P. Markowich, C. Ringhofer,  C. Schmeiser, \emph{Semiconductor equations}, Springer, Berlin, 1990.
%
\bibitem{Paredes}
A. Paredes, D.N. Olivieri, H. Michinel, \emph{From optics to dark matter: A review on nonlinear Schr\"odinger-Poisson systems} Physica D: Nonlinear Phenomena \textbf{403}, 132301, 2020.

\bibitem{RS}
C. Ringhofer, J. Soler, \emph{Discrete Schr\"odinger-Poisson systems preserving energy and mass}, Appl. Math. Lett. \textbf{13}, 27--32, 2000.
%
\bibitem{Shukla}
P.K. Shukla, B. Eliasson,  \emph{Colloquium: Nonlinear collective interactions in quantum plasmas with degenerate electron fluids}, Rev. Mod. Phys. \textbf{83}, 885, 2011.
\bibitem{Tha}
{\color{black}
M. Thalhammer, \emph{Convergence analysis of high-order time-splitting pseudo-spectral methods for nonlinear Schr\"odinger equations}, SIAM J.
Numer.  Anal. \textbf{50},  3231--3258, 2012.}

\bibitem{TM}
P. Tod, I.M. Moroz, \emph{An analytical approach to the Schr\"odinger-Newton equations}, Nonlinearity \textbf{12}, 201--216, 1999.
%
\bibitem{UKH}
C. Uhlemann, M. Kopp, and T. Haugg, \emph{Schr\"odinger method as N-body double and UV completion of dust}, Phys. Rev. D \textbf{90}, 023517, 2014. 
%
\bibitem{WK}
L. Widrow, N. Kaiser, \emph{Using the Schr\"odinger equation to simulate collisionless matter}, Astrophys. J. Lett. \textbf{416}, L71, 1993.
%
\bibitem{ZD}
Y. Zhang and X. Dong, \emph{On the computation of ground state and dynamics of Schr\"odinger-Poisson-Slater system}, J. Comput. Phys. \textbf{230}, 2660--2676, 2011.
%
\bibitem{Zhang}
Y. Zhang, \emph{Optimal error estimates of compact finite difference discretizations for the Schr\"odinger-Poisson system},  Commun. Commut. Phys. \textbf{13}, 1357--1388, 2015.
%
\bibitem{ZZM}
P. Zhang, Y. Zheng, N. Mauser, \emph{The limit from the Schr\"odinger-Poisson to the Vlasov-Poisson equations with general data in one dimension}, Comm. Pure Appl. Math.\textbf{55}, 582--632, 2002.
%
\bibitem{Zouraris1}
{\color{black} 
G. Zouraris, \emph{On the convergence of a linear two-step finite element method for the nonlinear
Schr\"odinger equation}, M2AN Math. Model. Numer. Anal. \textbf{35}, 389--405, 2001.}
%
\bibitem{Zouraris}
G.E. Zouraris, \emph{Error estimations pf the Besse relaxation scheme for a semilinear heat equation}, ESAIM: Math. Model. Numer. Anal. \textbf{55} 301--328, 2021.
%


\end{thebibliography}
\end{document}